\documentclass[11pt]{article}
\usepackage[margin=1in]{geometry}
\usepackage[T1]{fontenc}
\usepackage{lmodern,microtype}
\usepackage{amsmath,amssymb,amsthm}
\usepackage{needspace}
\usepackage{tikz}
\usepackage[hidelinks]{hyperref}

\newtheorem{theorem}{Theorem}[section]
\newtheorem{lemma}[theorem]{Lemma}
\newtheorem{corollary}[theorem]{Corollary}
\newtheorem{problem}{Problem}
\newcommand{\cut}{\partial}
\newcommand{\Ctwo}{\mathcal C_2}
\newcommand{\Htwo}{\mathcal H_2}

\title{Three-critical hypergraphs of maximum local edge connectivity two}
\author{Jun Yue\\[0.35em]
\small School of Mathematical Sciences, Tiangong University\\
\small Tianjin, China\\[0.2em]
\small\textit{Email address}:
\href{mailto:yuejun06@126.COM}{\texttt{yuejun06@126.COM}}}
\date{}

\begin{document}
\maketitle
\vspace{-2em}

\begin{abstract}
For a hypergraph \(G\), let \(\chi(G)\) denote its
chromatic number and \(\lambda(G)\) the maximum, over pairs of
distinct vertices, of the minimum size of an edge cut separating
the pair. A hypergraph is \(3\)-critical if its chromatic number
is three and every proper subhypergraph is two-colorable.
In this paper, we give a recursive characterization of the \(3\)-critical hypergraphs
with \(\lambda(G)\le2\). We prove that these are precisely
the hypergraphs generated from ordinary odd cycles and an explicit family
of labeled-tree hypergraphs by admissible Toft splits. This characterization
includes a necessary and sufficient condition for a prescribed
split to preserve \(\lambda\le2\). The condition involves only cuts in the
first factor and the two traces on the split edge. More
generally, we establish an exact cut formula for splits at
degree-\(k\) vertices of \((k+1)\)-critical hypergraphs,
for every \(k\ge2\). We further prove that the smallest member
of this critical class outside the family generated by
hyperwheels and both forms of the Haj\'os join has eight
vertices and is unique up to isomorphism. Together, these results
resolve two issues left open by Schweser, Stiebitz, and Toft in
\emph{Journal of Combinatorics} (2022):
the characterization of hypergraphs with \(\lambda(G)=2\)
and \(\chi(G)=3\), and whether hyperwheels and Haj\'os joins
suffice to generate all \(3\)-critical hypergraphs in this
connectivity range.
\end{abstract}

\section{Introduction}

All hypergraphs are finite and have distinct edges of size at
least two. For a hypergraph \(G\), write \(V(G)\) and \(E(G)\)
for its vertex set and edge set, and \(d_G(v)\) for the number
of edges containing \(v\). We omit the subscript when \(G\) is
clear. An edge of size two is \emph{ordinary}, and we write
\(uv\) for \(\{u,v\}\). A subhypergraph is obtained by retaining
a set of vertices together with a selection of the edges contained
in that vertex set. For \(U\subseteq V(G)\), the induced subhypergraph
\(G[U]\) has vertex set \(U\) and edge set
\(\{e\in E(G):e\subseteq U\}\). A coloring is \emph{proper}
if no edge is monochromatic, and \(\chi(G)\) is the least
number of colors in a proper coloring of \(G\). We call \(G\)
\emph{\(k\)-critical} if \(\chi(G)=k\) and \(\chi(H)<k\)
for every proper subhypergraph \(H\) of \(G\).

For \(U\subseteq V(G)\), the edge cut \(\cut_G(U)\) consists
of the edges meeting both \(U\) and \(V(G)\setminus U\).
The local edge connectivity of distinct vertices \(x,y\)
and the maximum local edge connectivity of \(G\) are,
respectively,
\[
 \lambda_G(x,y)=
 \min_{\substack{U\subseteq V(G)\\x\in U,\ y\notin U}}
 |\cut_G(U)|,
 \qquad
 \lambda(G)=\max_{\substack{x,y\in V(G)\\x\ne y}}
 \lambda_G(x,y).
\]
We put \(\lambda(G)=0\) if \(|V(G)|\le1\).
For \(v\in V(G)\), let \(G\div v\) be the hypergraph on
\(V(G)\setminus\{v\}\) whose edges are the distinct sets
\(e\setminus\{v\}\) of size at least two, for \(e\in E(G)\).
A vertex \(v\) of a connected hypergraph \(G\) is
\emph{separating} if \(G\div v\) is disconnected.
A \emph{block} of \(G\) is a maximal connected subhypergraph
with no separating vertex.

Toft \cite{toft-thesis} proved that
\(\chi(G)\le\lambda(G)+1\), an analogue of the chromatic
bound in terms of maximum degree. For ordinary graphs,
Aboulker et al.\ \cite{abht} and Stiebitz and Toft
\cite{st} studied coloring under local connectivity
constraints. Bastida and Brettell \cite{bb} considered
the corresponding list-coloring problem.
For hypergraphs with \(k=\lambda(G)\ge3\), the equality
case was determined by Schweser, Stiebitz, and Toft \cite[Theorem~2]{sst},
who showed that \(\chi(G)=k+1\) if and only if \(G\) contains
a block in \(\mathcal H_k\). The class \(\mathcal H_k\)
is generated by Haj\'os joins from the odd wheels when
\(k=3\) and from \(K_{k+1}\) when \(k\ge4\).

The case \(k=2\) leads to the class \(\Ctwo\) of
\(3\)-critical hypergraphs \(G\) satisfying \(\lambda(G)\le2\).
By Toft's bound, every member of \(\Ctwo\) has
\(\lambda(G)=2\). When the cut decomposition of Schweser,
Stiebitz, and Toft \cite[Theorem~12]{sst}
produces two critical factors from a member of \(\Ctwo\),
both factors again belong to \(\Ctwo\). In their concluding
remarks, the authors identify the inverse splitting operation as an
obstacle to completing the characterization
\cite[Section~6]{sst}. Schweser also
discusses this question in his thesis
\cite[Section~5.4, p.~104]{schweser-thesis}.

\begin{problem}\label{prob:splitting}
Characterize the hypergraphs \(G\) with \(\lambda(G)=2\)
and \(\chi(G)=3\).
\end{problem}

\phantomsection\label{subsec:recursive}
The inverse construction is given by Toft splitting
\cite{toft-critical,toft-hypergraphs,sst,book}.
Write \(G-e\) for the deletion of an edge \(e\) and
\(G-v\) for the deletion of a vertex \(v\) together with
all its incident edges. Let \(G_1,G_2\) be vertex-disjoint
hypergraphs. Choose
\(e\in E(G_1)\) and a vertex \(v\in V(G_2)\) with
\(d_{G_2}(v)=k\ge2\). Write the edges incident with \(v\) as
\(\{v\}\cup A_i\) for \(1\le i\le k\), where
\(A_i\subseteq V(G_2)\setminus\{v\}\). Choose nonempty sets
\(S_1,\ldots,S_k\subseteq e\) such that
\(\bigcup_{i=1}^k S_i=e\). We call the sets \(S_i\) the
\emph{traces} on \(e\), and allow them to overlap.
The split \(G^\ast\) has vertex set
\(V(G_1)\cup(V(G_2)\setminus\{v\})\) and edge set
\[
 E(G^\ast)
 =\bigl(E(G_1)\setminus\{e\}\bigr)
  \cup E(G_2-v)
  \cup\{S_i\cup A_i:1\le i\le k\}.
\]
If \(G_1\) and \(G_2\) are \((k+1)\)-critical, then \(G^\ast\)
is \((k+1)\)-critical \cite[Theorem~14]{sst}.
In particular, for two factors in \(\Ctwo\), the split
belongs to \(\Ctwo\) precisely when \(\lambda(G^\ast)\le2\).

A split can preserve criticality while violating the connectivity
bound; this can occur even for a split of two triangles. Let \(G_1\) be the triangle on \(\{a,b,c\}\),
let \(G_2\) be the triangle on \(\{v,x,y\}\), and take
\(e=ab\) and \(S_1=S_2=e\). The split inserts the edges
\(\{a,b,x\}\) and \(\{a,b,y\}\) and retains the path \(acb\).
Every cut separating \(a\) from \(b\) crosses both inserted
edges and at least one edge of this path. Hence
\(\lambda(G^\ast)\ge\lambda_{G^\ast}(a,b)\ge3\).

For \(S,U\subseteq V(G_1)\), the bipartition
\(U\mid V(G_1)\setminus U\) \emph{splits} \(S\) if both
parts meet \(S\). A Toft split at a vertex of degree
two is \emph{admissible} if every pair of distinct
vertices \(x,y\in V(G_1)\) admits a set \(U\subseteq V(G_1)\)
with \(x\in U\), \(y\notin U\), and \(|\cut_{G_1}(U)|=2\)
whose bipartition splits at most one of \(S_1,S_2\). The starting hypergraphs
for the recursion are the ordinary odd cycles and the
hypergraphs constructed from labeled trees in
Section~\ref{sec:bases}.

\begin{theorem}\label{thm:main}
A prescribed Toft split of two members of \(\Ctwo\) belongs
to \(\Ctwo\) if and only if it is admissible. Moreover,
\(\Ctwo\) is exactly the class generated from the ordinary
odd cycles and the tree constructions of
Section~\ref{sec:bases} by admissible splits. Consequently,
a connected hypergraph \(G\) with \(\lambda(G)\le2\)
satisfies \(\chi(G)=3\) if and only if it has a block in
\(\Ctwo\).
\end{theorem}

The block characterization answers
Problem~\ref{prob:splitting} componentwise. Indeed, the
chromatic number of a hypergraph is the maximum of the
chromatic numbers of its connected components, each of
which inherits the bound \(\lambda\le2\). A block of a
component is also a block of the whole hypergraph.
The role of the traces is made precise by
Lemma~\ref{lem:transfer}. To describe its two-color case, fix a
split of two members of \(\Ctwo\) and a bipartition \(U\mid V(G_1)\setminus U\) of the first factor.
The minimum size of an edge cut in \(G^\ast\) extending
this bipartition is \(|\cut_{G_1}(U)|\), plus one precisely
when the bipartition splits both traces. The same lemma
gives an exact cut formula for splits of two
\((k+1)\)-critical hypergraphs at a vertex of degree \(k\),
for every \(k\ge2\).
The decomposition into critical factors and the preservation of
criticality under a split are supplied by the results of Schweser,
Stiebitz, and Toft cited above. The recursive rule combines these
results with the cut criterion. Once the factors and traces are
specified, the criterion determines exactly whether the split also preserves
the required bound on maximum local edge connectivity.

The role of Haj\'os joins in the classes \(\mathcal H_k\)
for \(k\ge3\) suggests a corresponding generation rule
at \(k=2\). A \emph{hyperwheel} consists of a rim edge
and an apex joined by an ordinary edge to each rim vertex.
A rim of size two gives a triangle. The hyperwheels with
at least four vertices are precisely the hypergraphs
obtained from the one-node trees in Section~\ref{sec:bases}.
To form a Haj\'os join of vertex-disjoint hypergraphs
\(G_1,G_2\), choose \(e_i\in E(G_i)\) and \(v_i\in e_i\)
for \(i=1,2\). Delete \(e_1,e_2\), identify \(v_1,v_2\)
to a vertex \(z\), and add either \(e^\ast\) or
\(e^\ast\cup\{z\}\), where
\[
 e^\ast=(e_1\setminus\{v_1\})\cup(e_2\setminus\{v_2\})
\]
\cite{toft-thesis,toft-critical,sst,book}.
Let \(\Htwo\) be the class generated from hyperwheels
by both forms of the join.

Schweser, Stiebitz, and Toft
\cite[Section~6]{sst} observed that
\(\Htwo\subseteq\Ctwo\). They also exhibited examples
obtained by adding an edge on the leaves of certain
rooted trees. These examples also belong to \(\Htwo\), and
the authors left open whether the reverse inclusion holds.

\begin{problem}\label{prob:hajos}
Is there a hypergraph in \(\Ctwo\) that does not belong to
\(\Htwo\)?
\end{problem}

\begin{theorem}\label{thm:hajos}
The minimum order of a hypergraph in \(\Ctwo\setminus\Htwo\)
is eight. Up to isomorphism, there is a unique eight-vertex
member \(Q\) of \(\Ctwo\setminus\Htwo\). Its ordinary
edges are \(ua,ub,uc,ud,wp,wq\), and its other edges are
\(\{a,b,w\}\) and \(\{c,d,p,q\}\).
\end{theorem}

Figure~\ref{fig:obstruction} shows \(Q\). All starting
hypergraphs of the recursion belong to \(\Htwo\)
(Lemma~\ref{lem:hajos-bases}), and a Toft split of two
members of \(\Ctwo\) with a singleton trace is a Haj\'os
join of its factors (Lemma~\ref{lem:singleton}). It follows
that \(\Ctwo\) is generated from hyperwheels by Haj\'os
joins and admissible splits with \(|S_1|,|S_2|\ge2\), as
stated in Corollary~\ref{cor:two-traces}.

The remainder of the paper is organized as follows.
Section~\ref{sec:bases} records the elementary properties
of critical hypergraphs and characterizes the starting
hypergraphs. Section~\ref{sec:splits} establishes the cut
formula, the connectivity criterion, and the recursive
characterization. Section~\ref{sec:hajos} relates the
starting hypergraphs to Haj\'os joins, classifies the
small members, and determines the minimum order of an
obstruction together with uniqueness at that order.
Section~\ref{sec:related} compares the cut formula with
related results and proves that \(|E(G)|=|V(G)|\) for
every \(G\in\Ctwo\). It also establishes a sparsity
property of blocks with maximum local edge connectivity
at most two and relates the edge-count equality to
the theory of square critical hypergraphs.

\Needspace{10\baselineskip}
\section{The starting hypergraphs}\label{sec:bases}
\label{subsec:critical-properties}\label{subsec:labeled-trees}

We begin by collecting the properties of critical hypergraphs needed
in the structural arguments that follow. The degree and cut bounds control
the factors of a decomposition, while the edge-deletion criterion
will verify the criticality of the constructed hypergraphs.
The coloring bounds, in turn, identify the cases from which the
recursive construction starts.

\begin{lemma}[{\cite[Proposition~5(a)--(d)]{sst}}]
\label{lem:critical-properties}
Let \(k\ge1\), and let \(G\) be a \((k+1)\)-critical
hypergraph. Then \(G\) is connected, has no separating
vertex, and every vertex has degree at least \(k\).
Every nontrivial vertex partition has at least \(k\)
crossing edges. Moreover, no edge of \(G\) properly
contains another edge of \(G\).
\end{lemma}

The cut assertion is Toft's edge-connectivity bound, recorded
with the other properties in
\cite[Proposition~5(a)--(d)]{sst}.
In the applications below, it will provide a lower bound on
the size of every cut of a critical factor.

\begin{lemma}[{\cite[Proposition~4]{sst}}]
\label{lem:edge-deletion-input}
Let \(G\) be a connected hypergraph with at least one edge,
and let \(k\ge1\). Then \(G\) is \((k+1)\)-critical if
and only if \(\chi(G)>k\) and \(\chi(G-e)\le k\) for
every \(e\in E(G)\).
\end{lemma}

The connectedness hypothesis in
\cite[Proposition~4]{sst} is important when
applying this criterion; the same formulation appears in
\cite[Proposition~3.6]{sstorch}. For the constructions considered here, it is natural to
verify criticality by checking each edge deletion. We therefore record the
version for hypergraphs with no isolated vertices and explain why
connectedness follows from the remaining assumptions.

\begin{corollary}[{cf.~\cite[Proposition~4]{sst}}]\label{cor:edge-deletion-criterion}
Let \(G\) have no isolated vertices. Then \(G\) is
\(3\)-critical if and only if it is not two-colorable
and \(G-e\) is two-colorable for every \(e\in E(G)\).
\end{corollary}

\begin{proof}
The forward implication follows directly from criticality. To prove the
reverse implication, we first verify the connectedness hypothesis
of Lemma~\ref{lem:edge-deletion-input}. Since \(G\) is not
two-colorable, at least one connected component is not
two-colorable. If another component existed, it would
contain an edge because \(G\) has no isolated vertices.
Deleting an edge of that other component would leave the
non-two-colorable component unchanged, contradicting the
hypothesis on edge deletion. Hence \(G\) is connected.
It also has an edge, since an edgeless hypergraph is
two-colorable. Lemma~\ref{lem:edge-deletion-input} with
\(k=2\) now proves that \(G\) is \(3\)-critical.
\end{proof}

\begin{theorem}[{\cite{toft-thesis}}]\label{thm:toft-bound}
Every finite hypergraph \(G\) satisfies
\(\chi(G)\le\lambda(G)+1\).
\end{theorem}

Toft's bound \cite{toft-thesis} will give the upper bound on
the chromatic number in the block characterization. It is
also recorded in \cite[Section~1 and the beginning of Section~4]{sst}, where it connects critical subhypergraphs
with maximum local edge connectivity.

\emergencystretch=1em
\begin{theorem}[{\cite[Theorem~3]{ksw}}]
\label{thm:hypergraph-brooks}
Let \(G\) be a connected hypergraph, and put
\(\Delta(G)=\max_{x\in V(G)}d_G(x)\).
If \(G\) is neither a complete graph, nor an ordinary
odd cycle, nor a hypergraph consisting of a single edge,
then \(\chi(G)\le\Delta(G)\).
In particular, a \(3\)-critical hypergraph with maximum
degree at most two is an ordinary odd cycle.
\end{theorem}

The bound in \cite[Theorem~3]{ksw} supplies the first case
of the starting-hypergraph classification. Its critical form
is also contained in \cite[Lemma~6(d)]{sst};
to handle the remaining case, we must describe the incidences around
the unique vertex of higher degree.

The starting hypergraphs consist of the ordinary odd cycles
and a family of labeled-tree hypergraphs, which we now define
explicitly. Let \(T\) be a nonempty tree. Assign a distinct vertex \(p_{st}=p_{ts}\) to
each edge \(st\in E(T)\) and a set \(P_t\) of private vertices
to each node \(t\in V(T)\). The private sets are pairwise
disjoint and avoid all the vertices assigned to tree edges.
Writing \(N_T(t)\) for the set of neighbors of \(t\), put
\[
 B_t=P_t\cup\{p_{st}:s\in N_T(t)\},
 \qquad |B_t|\ge2\quad(t\in V(T)).
\]
Choose a set \(M\subseteq V(T)\) of marked nodes. Require a
proper labeling \(\sigma\colon V(T)\to\{0,1\}\) such that
\(\sigma(t)=0\) for \(t\in M\) and \(\sigma(t)=1\) whenever
\(P_t\ne\varnothing\). Here proper means that adjacent nodes
have different labels.

The vertex set consists of the assigned vertices, the private
vertices, and one new vertex \(v\), the \emph{apex}. For each
\(t\in V(T)\), include the \emph{node edge}
\[
 e_t=
 \begin{cases}
 B_t\cup\{v\}, & t\in M,\\
 B_t, & t\notin M.
 \end{cases}
\]
For each \(p\in\bigcup_{t\in V(T)}P_t\), also include the
\emph{spoke} \(\{v,p\}\). Require
\(d_G(v)=|M|+\sum_{t\in V(T)}|P_t|\ge3\).
The node edges are distinct because \(B_s\) and \(B_t\) share
at most one vertex for \(s\ne t\), and each has at least two
vertices. Every node edge contains at least two vertices other
than \(v\), so it is distinct from every spoke. Each \(p_{st}\)
belongs to \(e_s\) and \(e_t\), and each private vertex belongs
to its node edge and its spoke. Thus every vertex other than
\(v\) has degree two.

Theorem~\ref{thm:hypergraph-brooks} identifies the starting
hypergraphs whose vertices all have degree two. The one-apex
description is related to the low-degree structure in
\cite[Theorem~4]{ksw} and \cite[Lemmas~1.6 and 1.7]{ks}.
More precisely, Lemma~1.6(b), (e), and (f) of the latter paper
states that shrinking to the low vertices gives a Gallai forest,
and that distinct edges meeting the high vertices have distinct
low-vertex traces whenever those traces have size at least two;
these traces are bridges. A Gallai forest has only complete
graphs, ordinary odd cycles, and single edges as blocks.
Here a bridge is an edge whose vertices lie in distinct
components after its deletion.

Kostochka's \(Z_2\)-description
\cite[Lemma~4 and Property~8]{kost} gives further information
about the one-apex case. The converse in that description is stated only for
\(k\ge3\), so the corresponding two-color converse remains to be established
here. The direct argument below proves both the labeling
conditions and the criticality of the construction. It also
determines the possible colors at the boundary of each branch,
information used again in Section~\ref{sec:hajos}.
Theorem~9 of the same paper concerns the more restrictive class
of \(2\)-special hypergraphs, rather than all of \(\Ctwo\).
In that terminology, \(2\)-special means non-two-colorable with
a vertex of degree at most one in every proper nonempty
subhypergraph.

\begin{lemma}\label{lem:bases}
A \(3\)-critical hypergraph with no vertex of degree greater
than two is an ordinary odd cycle. A \(3\)-critical hypergraph
with exactly one vertex of degree greater than two is obtained
by the labeled-tree construction above. Conversely, all these
starting hypergraphs belong to \(\Ctwo\).
\end{lemma}

\begin{proof}
The first assertion is the critical case of
Theorem~\ref{thm:hypergraph-brooks}. For the second assertion,
suppose that \(v\) is the unique vertex of degree greater
than two. Lemma~\ref{lem:critical-properties} implies that
every other vertex has degree exactly two.
\emph{The auxiliary multigraph.}
Call the ordinary edges containing \(v\) spokes. To recover
the underlying tree and its private vertices from these incidences,
form a multigraph
\(D\) with one node \(t\) for each edge \(e_t\) of \(G\) other
than a spoke. Every vertex of degree two is in two distinct
edges and is in at most one spoke, since the spokes all contain
\(v\) and \(G\) has distinct edges. A degree-two vertex in two
node edges gives an edge of \(D\) between their nodes. A vertex
in a spoke is recorded as private at the node of its other
edge. Thus \(D\) has no loops, but it may have parallel edges.

Suppose \(D\) contains a cycle of length \(\ell\), where a pair
of parallel edges counts as a cycle of length two. In cyclic
order, write \(t_1,\ldots,t_\ell\) for its nodes and
\(p_1,\ldots,p_\ell\) for the vertices of \(G\) represented by
its successive edges, so \(p_i\in e_{t_i}\cap e_{t_{i+1}}\),
with indices taken modulo \(\ell\). Each \(p_i\) has degree two
and belongs only to those two edges. Delete all the edges
\(e_{t_i}\) from \(G\) and take a proper two-coloring of the
remaining hypergraph, keeping all its vertices. Recoloring
\(p_1,\ldots,p_\ell\) therefore affects none of the remaining
edges.
If \(\ell\) is even, alternate the colors on
\(p_1,\ldots,p_\ell\), so \(p_{i-1}\) and \(p_i\) have
different colors at every node \(t_i\). This gives a proper
two-coloring of \(G\). If \(\ell\) is odd and some node edge
has a vertex besides its two cycle vertices, choose the cyclic ordering so that this edge is
\(e_{t_1}\), and let \(x\) be a vertex of that edge
other than its two cycle vertices. The vertex \(x\) is outside
\(\{p_1,\ldots,p_\ell\}\), since each \(p_i\) lies only in its
two consecutive node edges. Alternate the colors along
\(p_1,\ldots,p_\ell\), starting with the color opposite to
that of \(x\). Then \(p_\ell\) and \(p_1\) have the same
color, opposite to that of \(x\), and each other node edge
contains two differently colored cycle vertices. This again
gives a proper two-coloring of \(G\).

The only remaining possibility is that \(\ell\) is odd and
\(e_{t_i}=\{p_{i-1},p_i\}\) for every \(i\). These edges form
an ordinary odd cycle, which must be all of \(G\) by
\(3\)-criticality. This is impossible because
\(d_G(v)>2\). The cycle argument therefore shows that \(D\)
is acyclic. Moreover, it has at least one node because every vertex other than \(v\) has degree
two and is in at most one spoke. If \(D\) had several components,
the node edges in each component, together with the spokes of
its private vertices, would form a proper subhypergraph.
These subhypergraphs cover all edges and share at most \(v\).
Choosing their two-colorings to agree at \(v\) would give a
proper two-coloring of \(G\). Hence \(D\) is a tree. Put \(T=D\).
Mark a node \(t\) if \(v\in e_t\) and let \(P_t\) be its
private set. Each node edge has at least two vertices other
than \(v\), because it has size at least two and is not a
spoke. This gives \(|B_t|\ge2\), with \(B_t=e_t\setminus\{v\}\)
as in the construction.

\emph{The labeling.}
With the underlying tree recovered, we determine the labels
from the colors permitted on its branches. For an oriented
tree edge \(st\), let \(G_{s\to t}\) have
the node edges indexed by the component of \(T-st\) containing
\(s\), together with the spokes of the private vertices at
those nodes. Its vertex set is the union of these edges and
\(\{v\}\), so it includes the boundary vertex \(p_{st}\).
Let \(L_{s\to t}\) be the set of colors attained by \(p_{st}\)
in proper two-colorings of \(G_{s\to t}\) with \(v\) colored
\(0\). Each branch is a proper subhypergraph of \(G\), since it
omits the node edge \(e_t\). Criticality therefore supplies
a proper two-coloring of the branch. If necessary, interchange
its two colors to give \(v\) color \(0\); this is also possible
when \(v\) is isolated in the branch. Consequently,
\[
 L_{s\to t}\ne\varnothing
 \qquad\text{for every oriented tree edge }st.
\]
The existence of branch colorings has therefore been established independently
of the restrictions on their possible boundary colors, which we derive next.

Now fix a node \(t\). The branches \(G_{s\to t}\) for
\(s\in N_T(t)\) share only \(v\). Each branch meets \(e_t\)
only in \(p_{st}\) and possibly \(v\). Consequently we can
choose a color from each \(L_{s\to t}\) independently.
Assign color \(1\) to the vertices of \(P_t\) to satisfy
their spokes. The resulting coloring satisfies every edge
except possibly \(e_t\). Since \(G\) is not two-colorable,
every such choice makes \(e_t\) monochromatic. Because
\(|B_t|\ge2\), varying the color of just one vertex
\(p_{st}\) while fixing the other colors would break
monochromaticity. Thus there is a color \(c_t\) such that
\(L_{s\to t}=\{c_t\}\) for all \(s\in N_T(t)\).
The private vertices and the apex, when present in \(e_t\),
also have color \(c_t\). In particular, \(c_t=0\) at marked
nodes and \(c_t=1\) wherever \(P_t\ne\varnothing\). For a
one-node tree the same reasoning applies with no branches.

To compare the labels at adjacent nodes \(s,t\), consider a proper two-coloring of
\(G_{s\to t}\) with \(v\) colored \(0\). For each
\(u\in N_T(s)\setminus\{t\}\), its restriction to
\(G_{u\to s}\) gives \(p_{us}\) color \(c_s\). The private
vertices at \(s\) have color \(1\), and the apex, if in
\(e_s\), has color \(0\). Thus every vertex in
\(e_s\setminus\{p_{st}\}\) has color \(c_s\). This set is
nonempty by \(|B_s|\ge2\). The edge \(e_s\) is proper only
if \(p_{st}\) has color \(1-c_s\). Thus the restriction
argument gives
\[
 L_{s\to t}\subseteq\{1-c_s\}.
\]
The nonemptiness already established makes this inclusion
an equality. The node argument at \(t\) also gives
\(L_{s\to t}=\{c_t\}\), so \(c_t=1-c_s\).
It follows that \(\sigma(t)=c_t\) defines a proper
two-color labeling of \(T\), with the required labels
at marked nodes and at nodes with private vertices.

\emph{Criticality of the construction.}
For the converse, let \(G\) be constructed from a labeled
tree \((T,\sigma)\). Define its branches and their color
sets as above. We first determine these color sets and then use them to
verify non-two-colorability and the effect of every edge deletion.
For every oriented tree edge \(st\), we claim that
\begin{equation}\label{eq:base-branch-colors}
 L_{s\to t}=\{1-\sigma(s)\}.
\end{equation}
Induct on the number of nodes on the \(s\)-side of \(T-st\).
For each \(u\in N_T(s)\setminus\{t\}\), the induction
hypothesis gives \(p_{us}\) color \(1-\sigma(u)=\sigma(s)\)
and supplies a branch coloring with this value. The vertices
of \(P_s\) must have color \(1=\sigma(s)\) if any are present,
and \(v\) has color \(0=\sigma(s)\) if \(s\) is marked.
When the \(s\)-side has just one node, the same assertions
hold with no child branches. In every case,
\(|B_s|\ge2\) ensures that
\(e_s\setminus\{p_{st}\}\) contains at least one vertex.
All its vertices have color \(\sigma(s)\) in any proper
coloring of the branch. Hence
\(L_{s\to t}\subseteq\{1-\sigma(s)\}\).
To prove the reverse inclusion, choose the child branch colorings
whose existence is provided by the induction hypothesis, give all private vertices at \(s\)
color \(1\), and give \(p_{st}\) color \(1-\sigma(s)\).
The child branches share only the apex, whose color has
been fixed to \(0\), and their boundary vertices are
distinct. These colorings therefore combine consistently.
Every child edge and every spoke is proper, while
\(e_s\) contains \(p_{st}\) and a vertex of the opposite
color. This constructs a proper coloring of the entire
branch and proves \(1-\sigma(s)\in L_{s\to t}\).
The two inclusions prove \eqref{eq:base-branch-colors},
including the induction base with no child branches.

Suppose \(G\) has a proper two-coloring and normalize it so
that \(v\) has color \(0\). At any node \(t\),
\eqref{eq:base-branch-colors} gives every boundary vertex
\(p_{st}\) color \(1-\sigma(s)=\sigma(t)\). The private
vertices at \(t\) and the apex if present in \(e_t\) have
the same color \(\sigma(t)\), so \(e_t\) is monochromatic.
Thus \(G\) is not two-colorable.

If \(e_t\) is deleted, combine the colorings of the branches
at \(t\) and assign color \(1\) to all private vertices at
\(t\). This colors every remaining edge properly. If instead
the spoke \(\{v,p\}\) with \(p\in P_t\) is deleted, use the
same branch colorings and give \(p\) color \(0\), leaving
the other private vertices at \(t\) colored \(1\). Here
\(\sigma(t)=1\), so \(t\) is unmarked and all vertices of
\(e_t\setminus\{p\}\) have color \(1\). Since \(|B_t|\ge2\),
at least one such vertex remains and \(e_t\) is proper.
The vertex \(p\) belongs only to \(e_t\) and the deleted
spoke, so this change preserves every other edge. The same
edge-deletion arguments also apply when the tree has only one node.
The construction has no isolated vertices, so
Corollary~\ref{cor:edge-deletion-criterion} gives
\(3\)-criticality. An ordinary odd
cycle is \(3\)-critical by Corollary~\ref{cor:edge-deletion-criterion}, since
deleting any edge gives a path.
In either family, every pair of distinct vertices contains a
vertex of degree two. The singleton cut at that vertex
separates the pair and has size two. Thus
\(\lambda(G)\le2\), proving membership in \(\Ctwo\).
\end{proof}

\section{Two-edge cuts and their inverse}\label{sec:splits}

\subsection{\texorpdfstring{Critical factors and the transfer formula}{Critical factors and the transfer formula}}\label{subsec:critical-factors}\label{subsec:cut-transfer}

We begin with the criticality of a split and the critical
factors obtained from a two-edge cut. The operation and
its traces are those defined in Section~\ref{subsec:recursive};
the following theorem applies to every degree allowed in
that definition.

\begin{theorem}[{\cite[Theorem~14]{sst}}]
\label{thm:toft-splitting}
Let \(k\ge2\), let \(G_1,G_2\) be vertex-disjoint
\((k+1)\)-critical hypergraphs, and choose
\(e\in E(G_1)\) and \(v\in V(G_2)\) with
\(d_{G_2}(v)=k\). Every Toft split of these factors
with nonempty traces covering \(e\) is
\((k+1)\)-critical. The inserted edges form a cut of
size \(k\) between the two factor sides.
\end{theorem}

For the converse direction, the traces must be recovered
from a cut together with the coloring restrictions on its
two sides. The next lemma uses the two-color specialization of
\cite[Theorem~12(a)--(d)]{sst}.
We will need its conclusions for a prescribed two-edge
cut, rather than only the existence of a cut with the
stated properties. The argument following the
statement explains why the conclusion holds for every prescribed cut;
the construction of the critical factors is exactly the one in the cited
theorem. The last assertion is the case \(k=2\) of
Theorem~\ref{thm:toft-splitting}, namely
\cite[Theorem~14]{sst}.

\begin{lemma}[{\cite[Theorem~12(a)--(d) and Theorem~14]{sst}}]\label{lem:toft}
Let \(G\) be \(3\)-critical, and let \(X\mid Y\) be a
nontrivial vertex partition with exactly two crossing edges
\(F_1,F_2\).
After interchanging \(X,Y\) if necessary, put
\(S_i=F_i\cap X\) and \(A_i=F_i\cap Y\) for \(i=1,2\).
Then every proper two-coloring of \(G[X]\) makes
\(S_1\cup S_2\) monochromatic, and every proper two-coloring
of \(G[Y]\) makes \(A_1,A_2\) monochromatic in opposite
colors. In particular, \(A_1\cap A_2=\varnothing\).

If \(|S_1\cup S_2|\ge2\), define \(G_1\) by adding the
edge \(e=S_1\cup S_2\) to \(G[X]\). Define \(G_2\) by adding
a new vertex \(v\) to \(G[Y]\) and the edges
\(\{v\}\cup A_1,\{v\}\cup A_2\). Both \(G_1,G_2\) are
\(3\)-critical, and \(d_{G_2}(v)=2\).

Conversely, every Toft split of two \(3\)-critical
hypergraphs at a degree-two vertex, with the traces specified
in Section~\ref{subsec:recursive}, is \(3\)-critical.
\end{lemma}

\begin{proof}[Justification for a prescribed cut]
Both induced sides are proper subhypergraphs of \(G\),
so each has a proper two-coloring. For a coloring of one
side, distinguish three possible trace patterns:
\(\mathsf F\), at least one trace is nonmonochromatic;
\(\mathsf S\), both traces are monochromatic in the same
color; and \(\mathsf O\), both traces are monochromatic
in opposite colors. If one side admitted pattern
\(\mathsf F\), combine that coloring with any proper
coloring of the other side. One crossing edge is already
proper, and interchanging the colors on the other side
can make the remaining crossing edge proper. This would
two-color \(G\), a contradiction. Thus every coloring
on either side has pattern \(\mathsf S\) or \(\mathsf O\).

Next, two side colorings with the same pattern would combine to give
a proper coloring of the whole hypergraph: after a possible interchange, each trace on
one side has the color opposite to its counterpart on
the other side. Fix a coloring on one side. Every
coloring on the other must therefore have the opposite
pattern. Fixing one of those colorings in turn shows
that all colorings on the first side have the original
pattern. Consequently, after interchanging the sides of the prescribed
partition if necessary, every coloring on \(X\) has pattern
\(\mathsf S\) and every coloring on \(Y\) has pattern
\(\mathsf O\). This also gives
\(A_1\cap A_2=\varnothing\), since a vertex in the
intersection could not receive the two opposite colors.

These are precisely the coloring properties used in
the critical-factor construction in
\cite[Theorem~12(c),(d) and their proofs]{sst}.
That construction applies to the partition just
oriented and gives the stated factors, without replacing
the prescribed cut. The converse assertion is
Theorem~\ref{thm:toft-splitting} with \(k=2\).
\end{proof}

The trace-color conditions distinguish the two factors:
one acquires a new edge, while the other acquires a new
vertex. To control a prescribed inverse split, we must now
determine how cuts extend across these two sides and how
their local edge connectivities change.

The cut calculation uses uncrossing for the symmetric
submodular edge-cut function, a framework also used in
hypergraph minimum-cut algorithms \cite{cx}. When contracting a set of vertices,
we retain any parallel edges and count them with their
multiplicities, while discarding edges reduced to a single vertex. These conventions ensure that corresponding
cuts before and after contraction have the same number
of crossing edges. With these conventions for measuring cut size,
the following submodularity property applies equally to the
hypergraphs obtained by contraction.

\begin{lemma}[{\cite[Section~2.1, p.~2125]{cx}}]
\label{lem:cut-function-input}
Let \(G\) be a finite hypergraph, allowing parallel edges
counted with their multiplicities, and put
\(f(A)=|\cut_G(A)|\) for \(A\subseteq V(G)\).
Then \(f(A)=f(V(G)\setminus A)\), and for all
\(A,B\subseteq V(G)\),
\[
 f(A)+f(B)\ge f(A\cap B)+f(A\cup B).
\]
\end{lemma}

The weighted version of this property in \cite[Section~2.1, p.~2125]{cx}
also covers parallel edges: their multiplicities are added
as an edge weight. We use submodularity to extend a fixed
partition on one side of a global minimum cut without
splitting the other side. The next lemma formulates the precise extension property
that will be used in the transfer argument.

\begin{lemma}\label{lem:contraction}
Let \(X\mid Y\) be a global minimum cut of a hypergraph \(G\),
with \(X,Y\ne\varnothing\). For every \(U\subseteq X\),
there is a minimum-size cut extending \(U\mid X\setminus U\)
that puts all of \(Y\) on one side. Consequently, contracting
\(Y\) to one vertex preserves the local edge connectivity
of every pair of distinct vertices in \(X\), with parallel
edges counted with multiplicity.
\end{lemma}

\begin{proof}
Put \(f(Z)=|\cut_G(Z)|\) for \(Z\subseteq V(G)\).
By Lemma~\ref{lem:cut-function-input}, this function is
symmetric and submodular.

Choose \(Z\subseteq V(G)\) of minimum \(f(Z)\) subject
to \(Z\cap X=U\). If \(Z\) does not split \(Y\), this
minimum extension already has the required form.
Otherwise \(Z\cap Y\) and \(Y\setminus Z\)
are nonempty. Hence \(Z\cup X\) is a nonempty proper subset
of \(V(G)\), and global minimality gives
\(f(Z\cup X)\ge f(X)\). Submodularity yields
\[
 f(Z)+f(X)\ge f(Z\cap X)+f(Z\cup X)
             \ge f(U)+f(X).
\]
Thus \(f(U)\le f(Z)\). The set \(U\), viewed as a subset
of \(V(G)\), is an admissible extension with all of \(Y\)
in its complement. It therefore also attains the minimum.
A cut of the contracted hypergraph corresponds exactly to
a cut of \(G\) with \(Y\) wholly on one side, and each original
edge contributes the same crossing indicator to these two
cuts. Keeping multiplicities makes their sizes equal.
For a pair \(x,x'\in X\), the first assertion shows that
a minimum cut separating them can be chosen with \(Y\)
intact. Taking the minimum over these corresponding cuts therefore proves
the assertion about local edge connectivity.
\end{proof}

For a Toft split, Theorem~\ref{thm:toft-splitting} gives
criticality, and Lemma~\ref{lem:critical-properties} supplies
the lower bound on nontrivial cuts. Thus the cut between
the factor sides is a global minimum cut.
Lemma~\ref{lem:contraction} then reduces the extension problem
to the two possible placements of the opposite side. Comparing
their contributions gives the exact formula below.

\begin{lemma}\label{lem:transfer}
Let \(k\ge2\). Take vertex-disjoint \((k+1)\)-critical
hypergraphs \(G_1,G_2\), an edge \(e\in E(G_1)\), and
a vertex \(v\in V(G_2)\) of degree \(k\), whose incident
edges are \(\{v\}\cup A_1,\ldots,\{v\}\cup A_k\).
Let \(G\) be their Toft split with nonempty traces
\(S_1,\ldots,S_k\subseteq e\) covering \(e\). Put
\(X=V(G_1)\) and \(Y=V(G_2)\setminus\{v\}\). For
\(U\subseteq X\), define
\[
 p(U)=|\{i:S_i\subseteq U\}|,\qquad
 q(U)=|\{i:S_i\subseteq X\setminus U\}|.
\]
Then
\begin{equation}\label{eq:transfer}
 \min_{\substack{Z\subseteq V(G)\\Z\cap X=U}}
       |\cut_G(Z)|
 =|\cut_{G_1-e}(U)|+k-\max\{p(U),q(U)\}.
\end{equation}
Moreover,
\[
 \lambda_{G_1}(x,x')\le\lambda_G(x,x')
 \quad(x,x'\in X,\ x\ne x'),
 \qquad
 \lambda_G(y,y')=\lambda_{G_2}(y,y')
 \quad(y,y'\in Y,\ y\ne y'),
\]
and \(\lambda_G(x,y)=k\) for every \(x\in X\), \(y\in Y\).
\end{lemma}

\begin{proof}
The sets \(A_i\) are nonempty and pairwise distinct:
each incident edge has size at least two, and distinct
edges at \(v\) give distinct sets after its removal.
Thus the inserted edges \(S_i\cup A_i\) are distinct
and meet both sides \(X,Y\).
They are precisely the crossing edges of \(X\mid Y\), so
this cut has size \(k\). Theorem~\ref{thm:toft-splitting}
implies that \(G\) is \((k+1)\)-critical.
Lemma~\ref{lem:critical-properties} therefore shows that
every nontrivial cut has size at least \(k\), so
\(X\mid Y\) is a global minimum cut.
Fix \(U\subseteq X\). Lemma~\ref{lem:contraction} shows
that a minimum-size extension may put all of \(Y\) on one
side. Edges internal to \(Y\) then contribute zero, and
edges internal to \(X\) contribute
\(|\cut_{G_1-e}(U)|\). If \(Y\subseteq Z\), the inserted
edge \(S_i\cup A_i\) crosses precisely when \(S_i\not\subseteq U\),
giving contribution \(k-p(U)\). If \(Y\cap Z=\varnothing\),
it crosses precisely when \(S_i\not\subseteq X\setminus U\),
giving contribution \(k-q(U)\). Comparing these two values and taking the smaller one
proves \eqref{eq:transfer}. This also covers \(U=\varnothing\)
and \(U=X\), for which the minimum is zero.

If \(U\) does not split \(e\), all \(k\) traces lie on one
side, and the right side of \eqref{eq:transfer} equals
\(|\cut_{G_1}(U)|\). If \(U\) splits \(e\), the covering
condition \(\bigcup_iS_i=e\) implies
\(\max\{p(U),q(U)\}\le k-1\). The new contribution is then
at least one, which is the contribution of \(e\) to
\(\cut_{G_1}(U)\). Thus the minimum extension size is
always at least \(|\cut_{G_1}(U)|\). Taking the minimum
over all \(U\) separating \(x,x'\) proves the first
connectivity inequality.
For the second factor, apply Lemma~\ref{lem:contraction} with \(X,Y\) interchanged
to a pair \(y,y'\in Y\). Contracting \(X\) to a vertex
identified with \(v\) sends each inserted edge to
\(\{v\}\cup A_i\), discards the edges internal to \(X\),
and retains exactly the edges of \(G_2-v\). The contracted
hypergraph is therefore \(G_2\), proving the asserted
equality on \(Y\). Finally, the cut \(X\mid Y\) separates
every \(x\in X\) from every \(y\in Y\) and has size \(k\).
The lower bound on all nontrivial cuts gives the reverse
inequality, so \(\lambda_G(x,y)=k\).
\end{proof}

\subsection{Preserving maximum local edge connectivity}\label{subsec:connectivity-test}

The transfer formula determines the local edge connectivities for pairs
within the second factor and for pairs on opposite sides. For a pair in the first
factor, we must determine whether a minimum cut extends
without increasing its size. The following criterion expresses
this requirement entirely through cuts of the first factor
and the prescribed traces.

\begin{corollary}\label{cor:split-connectivity}
In the setting of Lemma~\ref{lem:transfer}, suppose
\(\lambda(G_1),\lambda(G_2)\le k\). Then \(\lambda(G)\le k\)
if and only if every pair of distinct vertices in \(X\)
is separated by a cut of size \(k\) in \(G_1\) with at least
\(k-1\) traces wholly on the same side.
\end{corollary}

\begin{proof}
For \(U\subseteq X\), let \(\delta_e(U)\) be \(1\) if \(U\)
splits \(e\), and \(0\) otherwise. Subtracting the contribution
of \(e\) in \eqref{eq:transfer} expresses the minimum
extension size as
\[
 |\cut_{G_1}(U)|+\varepsilon(U),\qquad
 \varepsilon(U)=k-\max\{p(U),q(U)\}-\delta_e(U)\ge0.
\]
If \(U\) does not split \(e\), then
\(\max\{p(U),q(U)\}=k\) and \(\varepsilon(U)=0\).
If \(U\) splits \(e\), then this maximum is at most \(k-1\),
and \(\varepsilon(U)=0\) precisely when it is \(k-1\).
Thus \(\varepsilon(U)=0\) is equivalent to at least \(k-1\)
traces lying wholly on the same side.
Every nontrivial cut of \(G_1\) has size at least \(k\).
Consequently, a pair in \(X\) has local edge connectivity
at most \(k\) in \(G\) precisely when it admits a separating
set \(U\) with \(|\cut_{G_1}(U)|=k\) and
\(\varepsilon(U)=0\). Pairs in \(Y\) have connectivity at
most \(k\) by Lemma~\ref{lem:transfer} and the hypothesis
on \(G_2\), and cross pairs have connectivity exactly \(k\).
These observations prove both directions of the stated criterion.
\end{proof}

\begin{corollary}\label{cor:admissibility-test}
A Toft split of two members of \(\Ctwo\) has maximum local
edge connectivity at most two if and only if it is
admissible. It suffices to check the admissibility condition
for pairs contained in
\[
 W=\{x\in V(G_1):d_{G_1}(x)\ge3\}\cup(S_1\cap S_2).
\]
Every admissible split satisfies \(|S_1\cap S_2|\le1\).
\end{corollary}

\begin{proof}
For \(k=2\), the minimum extension size in
\eqref{eq:transfer} is
\begin{equation}\label{eq:two-trace-transfer}
 \min_{\substack{Z\subseteq V(G)\\Z\cap X=U}}
       |\cut_G(Z)|
 =|\cut_{G_1}(U)|+\varepsilon(U),
\end{equation}
where \(\varepsilon(U)=1\) if \(U\mid X\setminus U\)
splits both \(S_1,S_2\), and \(\varepsilon(U)=0\) otherwise.
Indeed, if \(e\) is intact, then \(p(U)+q(U)=2\) and one
of these values is two. If \(e\) is split, at most one trace
can lie wholly on either fixed side. The value of
\(\max\{p(U),q(U)\}\) is zero exactly when both traces
are split, and it is one otherwise. Equation
\eqref{eq:two-trace-transfer}, or
Corollary~\ref{cor:split-connectivity}, gives the
admissibility equivalence.
Every vertex \(x\notin W\) of the first factor has
degree two, since \(G_1\) is \(3\)-critical. Its singleton cut has size two and
splits at most one trace, because \(x\notin S_1\cap S_2\).
Thus every pair with an endpoint outside \(W\) already
has a suitable separating cut. Finally, if distinct
vertices belong to \(S_1\cap S_2\), every cut separating
them splits both traces. Such a pair violates admissibility,
proving the intersection bound.
\end{proof}

The example of two triangles in the introduction extends
to every hyperwheel. Let \(e\) be its rim and take
\(S_1=S_2=e\), with the second factor a triangle. Every cut
separating two rim vertices crosses both inserted edges and
at least one spoke. Hence this split has maximum local edge
connectivity at least three.

\subsection{Blocks and the recursive characterization}\label{subsec:recursive-characterization}

To pass from the description of critical hypergraphs to the
block characterization, we use the attachment argument in
\cite[Proof of Theorem~2]{sst} in the range
\(\lambda\le2\). The argument compares a lower bound on cuts
inside the critical subhypergraph with the matching upper
bound in the ambient hypergraph. The next lemma states
these assumptions directly. The proof verifies both the attachment
property and the maximality assertion, making explicit why this argument
also applies when \(k=2\).

\begin{lemma}[{cf.~\cite[Proof of Theorem~2]{sst}}]\label{lem:embedded-block}
Let \(G\) be a hypergraph with \(\lambda(G)\le2\).
Suppose that \(H\subseteq G\) is a connected subhypergraph
on at least two vertices, has no separating vertex, and
every nontrivial edge cut of \(H\) has size at least two.
Then \(H\) is a block of \(G\).
\end{lemma}

\begin{proof}
We first recall that local edge connectivity is monotone under taking
subhypergraphs. If \(a,b\in V(J)\) for a subhypergraph
\(J\subseteq G\), then every cut of \(G\) separating them
restricts to a cut of \(J\) separating them, with no more
crossing edges. Taking a minimum cut of \(G\) gives
\(\lambda_J(a,b)\le\lambda_G(a,b)\), and hence
\(\lambda(J)\le2\).
Let \(J\) be a connected subhypergraph of \(G\) containing
\(H\), and delete \(E(H)\) from \(J\), retaining all vertices.
Every component \(C\) of the resulting hypergraph meets
\(V(H)\) in at most one vertex. Otherwise take distinct
\(a,b\in C\cap V(H)\). Every cut of \(J\) separating them
crosses at least two edges of \(H\), by the hypothesis
on \(H\). It also crosses an edge outside \(H\), since
\(a,b\) lie in the same connected component after all
edges of \(H\) are deleted. These are distinct crossing
edges, so \(\lambda_J(a,b)\ge3\), a contradiction.

Suppose \(J\) properly contains \(H\). It must have an edge
outside \(H\), since a connected enlargement of \(H\) cannot
add only vertices. Let \(C\) be the component containing
such an edge after the deletion of \(E(H)\). This component
meets \(V(H)\). Indeed, if it did not, no edge of \(H\)
would meet \(C\), and no edge outside \(H\) could join
\(C\) to its complement by the definition of a component.
This would disconnect \(J\). Thus \(C\cap V(H)=\{a\}\)
for some vertex \(a\).
Both \(C\setminus\{a\}\) and \(V(H)\setminus\{a\}\) are
nonempty. The first assertion follows because \(C\)
contains an edge of size at least two. Every edge of \(J\)
outside \(H\) lies in one component after the deletion
of \(E(H)\), and every edge of \(H\) avoids
\(C\setminus\{a\}\). After removing \(a\) from all edges,
no resulting edge meets both \(C\setminus\{a\}\) and
\(V(J)\setminus C\). These are nonempty parts of
\(V(J)\setminus\{a\}\), so \(J\div a\) is disconnected.
Thus every proper connected enlargement of \(H\) has a
separating vertex. Together with the hypotheses on this subhypergraph, this establishes precisely the
maximality required for \(H\) to be a block of \(G\).
\end{proof}

We now combine the critical-factor decomposition with
the connectivity criterion to prove the recursive description.
The block lemma then allows us to pass from this description
of critical hypergraphs to the characterization of the ambient hypergraph.

\begin{proof}[Proof of Theorem~\ref{thm:main}]
By Lemma~\ref{lem:toft}, every Toft split of two members of
\(\Ctwo\) is \(3\)-critical. By
Corollary~\ref{cor:admissibility-test}, its maximum local
edge connectivity is at most two precisely when the split
is admissible. This proves the first assertion. All starting
hypergraphs belong to \(\Ctwo\) by Lemma~\ref{lem:bases}.
Induction on the number of splits therefore shows that
every hypergraph generated by the stated operations
belongs to \(\Ctwo\).

For the reverse inclusion, induct on \(|V(G)|\) for
\(G\in\Ctwo\). If at most one vertex has degree greater
than two, Lemma~\ref{lem:bases} identifies \(G\) as a
starting hypergraph. Otherwise choose vertices \(a,b\)
of degree at least three. Every nontrivial cut has size
at least two, while \(\lambda_G(a,b)\le2\). There is
therefore a cut \(X\mid Y\) of size two separating \(a,b\).
Neither side is a singleton, since a singleton side
would consist of \(a\) or \(b\) and its cut size would
equal that vertex's degree.
Each side contains at least two vertices incident with
crossing edges. To see this, suppose all crossing edges
meet \(X\) in the same vertex \(x\). Since \(|X|\ge2\)
and \(Y\ne\varnothing\), the sets \(X\setminus\{x\}\)
and \(Y\) are nonempty. On removing \(x\) from all edges,
each former crossing edge lies wholly in \(Y\), and
each internal edge remains within its original side.
Thus \(G\div x\) is disconnected, contrary to the absence
of separating vertices in a \(3\)-critical hypergraph.
The same argument applies to \(Y\).

Orient the cut according to Lemma~\ref{lem:toft}. Its
boundary on \(X\) has size at least two, so the lemma
gives \(3\)-critical factors \(G_1,G_2\) with
\[
 |V(G_1)|=|X|\le |V(G)|-2,\qquad
 |V(G_2)|=|Y|+1\le |V(G)|-1.
\]
By construction, the corresponding Toft split is exactly
\(G\). Lemma~\ref{lem:transfer} gives
\(\lambda_{G_1}(x,x')\le\lambda_G(x,x')\le2\)
for distinct \(x,x'\in X\), and
\(\lambda_{G_2}(y,y')=\lambda_G(y,y')\le2\)
for distinct \(y,y'\in Y\). The new vertex \(v\) of
\(G_2\) has degree two, so its singleton cut separates
it from every vertex of \(Y\) with two crossing edges.
Thus both factors belong to \(\Ctwo\). Since the reconstructed
split satisfies \(\lambda(G)\le2\),
Corollary~\ref{cor:admissibility-test} shows that it
is admissible.
Both factors have smaller order, so the induction applies
to each and supplies a decomposition of \(G\).

Finally, let \(G\) be connected with \(\lambda(G)\le2\).
Theorem~\ref{thm:toft-bound} gives \(\chi(G)\le3\). If \(\chi(G)=3\),
choose a subhypergraph \(H\) minimal under inclusion
among those that are not two-colorable. Since
\(\chi(H)\le\chi(G)\le3\), this \(H\) is \(3\)-critical.
Monotonicity gives \(\lambda(H)\le2\), so \(H\in\Ctwo\).
Lemma~\ref{lem:critical-properties} gives connectedness,
the absence of a separating vertex, and a lower bound
of two on its nontrivial cut sizes.
Lemma~\ref{lem:embedded-block} therefore shows that
\(H\) is a block of \(G\). Conversely, if \(G\) has a block
in \(\Ctwo\), every proper coloring of \(G\) restricts
to a coloring of that block, so \(\chi(G)\ge3\).
Together with the upper bound, this gives \(\chi(G)=3\).
\end{proof}

\section{The smallest obstruction to Haj\'os generation}\label{sec:hajos}

The obstruction analysis uses the criticality and
local-connectivity properties of Haj\'os joins. We first
record these properties for both forms of the operation
defined in the introduction, and then determine the conditions under which
such a join can be reversed into critical factors.

\begin{lemma}[{\cite[Theorem~9(a), proof of Theorem~13, Claim~2, and Section~6]{sst}}]
\label{lem:hajos-input}
Let \(k\ge2\). A Haj\'os join of two \((k+1)\)-critical
hypergraphs is \((k+1)\)-critical. If both factors belong
to \(\Ctwo\), then their join also belongs to \(\Ctwo\).
Consequently, \(\Htwo\subseteq\Ctwo\).
\end{lemma}

The criticality assertion is
\cite[Theorem~9(a)]{sst}.
The forward implication in the proof of
\cite[Theorem~13, Claim~2]{sst} gives the
local-connectivity bound
\(\lambda(G_1\mathbin{\Delta}G_2)
\le\max\{\lambda(G_1),\lambda(G_2)\}\).
That part of the argument uses neither \(k\ge3\) nor
the reverse-criticality assertion; it applies to both
forms of the join at \(k=2\). Together with the criticality assertion, this bound proves
the closure assertion stated in the lemma.
The inclusion \(\Htwo\subseteq\Ctwo\) is also stated
explicitly in \cite[Section~6]{sst}.
For the inverse operation, we must also determine the
color relation between the removed edge and the identified
vertex. This relation is the subject of the next subsection.

\subsection{\texorpdfstring{Inverse joins and an eight-vertex obstruction}{Inverse joins and an eight-vertex obstruction}}

Call a pair \((z,e)\), where \(z\in V(G)\) and \(e\in E(G)\),
a \emph{mixed separator} if \((G-e)\div z\) is disconnected.
The mixed separator theorem
\cite[Theorem~10]{sst} concerns
\((k+1)\)-critical hypergraphs with \(k\ge3\). For
\(3\)-critical hypergraphs, a mixed separator must also
satisfy a coloring condition in order to recover critical
factors. The same restriction \(k\ge3\) also occurs in
the converse criticality assertion for Haj\'os joins
\cite[Theorem~9(b)]{sst}, so neither
result applies directly with \(k=2\).
Lemma~\ref{lem:inverse} identifies the necessary and sufficient
color relation in this case. Its proof also verifies the
criticality, smaller order, and local-connectivity bound
of the recovered factors.

\begin{lemma}\label{lem:inverse}
Let \(G\) be \(3\)-critical and let \((z,e)\) be a mixed
separator. There is a constant
\(\varepsilon_{z,e}\in\{0,1\}\), depending only on
\((G,z,e)\), such that every proper two-coloring
\(\varphi\colon V(G)\to\{0,1\}\) of \(G-e\) satisfies
\begin{equation}\label{eq:inverse-uniform-color}
 \varphi(x)=\varphi(z)\oplus\varepsilon_{z,e}
 \qquad\text{for every }x\in e,
\end{equation}
where \(\oplus\) denotes addition modulo two.
In particular, \(e\) is monochromatic in every such
coloring. The constant is independent both of the
coloring and of any partition of the components of
\((G-e)\div z\).

The hypergraph \(G\) is a Haj\'os join of two \(3\)-critical
hypergraphs with identified vertex \(z\) and new edge \(e\)
if and only if \(\varepsilon_{z,e}=0\), that is, if and
only if the common color on \(e\) equals the color of
\(z\) in every proper two-coloring of \(G-e\).
When this equality holds, every partition of the components
of \((G-e)\div z\) into two nonempty groups yields such
factors. Both factors have fewer vertices than \(G\) and
maximum local edge connectivity at most \(\lambda(G)\).
\end{lemma}

\begin{proof}
Every nontrivial cut of \(G\) has size at least two, so
\(G-e\) is connected. Otherwise a component of \(G-e\)
would define a nontrivial cut of \(G\) with at most the
single crossing edge \(e\). Every proper two-coloring
of \(G-e\) makes \(e\) monochromatic, since a coloring
making it nonmonochromatic would color \(G\).
Every component of \((G-e)\div z\) meets
\(e\setminus\{z\}\). If a component did not, restoring
\(e\) could not join it to another component after
shrinking \(z\), and \(G\div z\) would still be
disconnected. This contradicts the absence of separating
vertices in \(G\).

Fix a proper two-coloring \(\varphi_0\) of \(G-e\),
normalized by \(\varphi_0(z)=0\); it exists by criticality.
Choose \(x_0\in e\), and put
\(\varepsilon_{z,e}=\varphi_0(x_0)\).
Monochromaticity of \(e\) makes this value independent
of the choice of \(x_0\). We now show that this same value
determines the color relation for every coloring and every grouping
of the components.
Partition the components arbitrarily into two nonempty
groups, with vertex unions \(P_1,P_2\). For \(i=1,2\), put
\[
 J_i=(G-e)[P_i\cup\{z\}],\qquad
 D_i=(e\cap P_i)\cup\{z\},\qquad
 H_i=J_i+D_i.
\]
Both \(e\cap P_i\) are nonempty, so \(|D_i|\ge2\).
Every edge of \(G-e\) is contained in one
\(P_i\cup\{z\}\), because its vertices other than \(z\)
belong to one component after shrinking \(z\).
Moreover, each \(J_i\) is connected. Its components
away from \(z\) are the chosen components of
\((G-e)\div z\), and connectedness of \(G-e\) forces
each of them to be joined to \(z\).
The hypergraphs \(J_1,J_2\) admit proper two-colorings
by restriction from \(G-e\). Once \(z\) has color \(0\),
their colorings can be chosen independently and combined
to color \(G-e\). If \(z\in e\), monochromaticity forces
every vertex of \(e\) to have color \(0\). If \(z\notin e\),
fix a coloring of \(J_2\). Since \(e\cap P_2\ne\varnothing\)
and every combined coloring makes \(e\) monochromatic,
this fixed coloring forces a single color on
\(e\cap P_1\) in every coloring of \(J_1\).
Fixing one such coloring of \(J_1\) then forces the same
color on \(e\cap P_2\) in every coloring of \(J_2\).
Thus in either case there is a fixed \(c\in\{0,1\}\)
such that every coloring of either \(J_i\) with \(z\)
colored \(0\) gives \(e\cap P_i\) color \(c\).
The restrictions of the fixed coloring \(\varphi_0\)
are among these colorings of \(J_1,J_2\), so
\(c=\varepsilon_{z,e}\). Every coloring of \(G-e\)
with \(z\) colored \(0\) restricts to colorings of
these two hypergraphs, and hence gives all of \(e\)
this same color. Interchanging the two colors proves
\eqref{eq:inverse-uniform-color} for either color at
\(z\). Since \(\varphi_0\) was fixed before the
arbitrary partition was chosen, the value of
\(\varepsilon_{z,e}\) is independent of that partition.
It follows that the constructions below apply to every
partition of the components into two nonempty groups.

Consider first the case \(c=1\). Each coloring of \(J_i\) with \(z\) colored
\(0\) makes \(D_i\) proper, so \(H_i\) is two-colorable.
Every inverse Haj\'os representation with identified
vertex \(z\) and new edge \(e\) arises from a partition
of the components into two groups of this kind. Indeed,
after deleting the new edge, no edge joins the two factor
sides away from \(z\), so each component after shrinking
\(z\) belongs entirely to one factor. The deleted edge
of that factor must then be \(D_i\).
Consequently, \(c=1\) permits no representation with
two \(3\)-critical factors.

Now suppose \(c=0\). Every proper two-coloring of \(J_i\)
makes \(D_i\) monochromatic. In particular, \(D_i\)
is not already an edge of \(J_i\), and \(H_i\) is not
two-colorable. Deleting \(D_i\) leaves the two-colorable
hypergraph \(J_i\). Let \(a\in E(J_i)\), and take a
proper two-coloring of \(G-a\), normalized so that \(z\)
has color \(0\). Its restriction to \(J_{3-i}\) gives
\(e\cap P_{3-i}\) color \(0\). Since \(e\) is proper
in this coloring, some vertex of \(e\cap P_i\) has
color \(1\). Thus restriction to \(P_i\cup\{z\}\)
properly colors \(H_i-a\), including \(D_i\).
The connected hypergraph \(J_i\) has at least two
vertices, so \(H_i\) has no isolated vertices.
Corollary~\ref{cor:edge-deletion-criterion}
now gives \(3\)-criticality.
Take disjoint copies of \(H_1,H_2\), identify their copies
of \(z\), delete \(D_1,D_2\), and add their union with
\(z\) omitted or included according as \(z\notin e\) or
\(z\in e\). The resulting edge is exactly \(e\), and
all other edges are those of \(G-e\). This restores \(G\).
Also \(|V(H_i)|=|P_i|+1<|V(G)|\), since both groups
are nonempty.

It remains to verify the connectivity bound. Fix two vertices of \(H_i\)
and consider any cut of \(G\) separating them. Restrict this cut
to \(P_i\cup\{z\}\). Each retained edge of \(J_i\) has the same crossing
status in the restricted cut as in the original cut,
so its contribution is unchanged. If \(D_i\) crosses, its contribution
can be accounted for by an edge of \(G\) outside \(J_i\).
When \(e\) crosses, use \(e\) itself. Otherwise the
crossing of \(D_i\) places all vertices of \(e\) on
the side opposite to \(z\). Choose a vertex of
\(e\cap P_{3-i}\). Since \(J_{3-i}\) is connected,
some edge of \(J_{3-i}\) crosses between this vertex
and \(z\). Thus the restricted cut of \(H_i\) has
no more crossing edges than the original cut of \(G\).
Taking a minimum cut for each pair proves
\(\lambda(H_i)\le\lambda(G)\).
\end{proof}

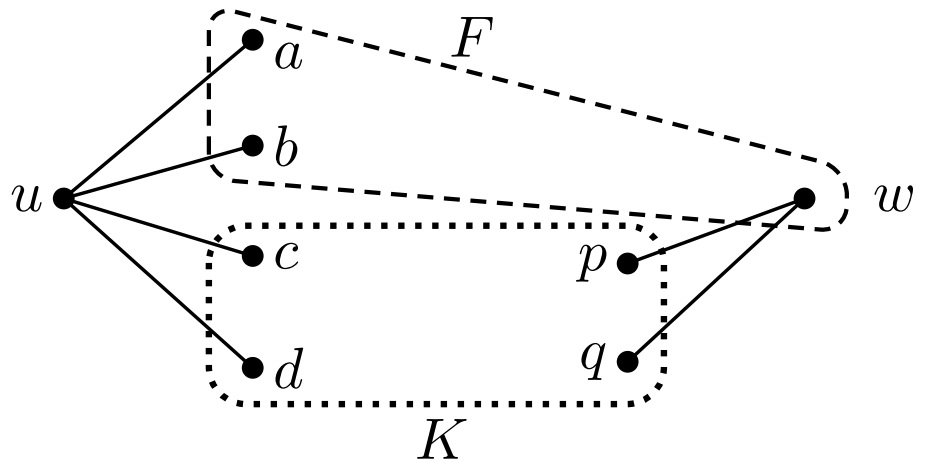
\begin{figure}[htbp]
\centering
\begin{tikzpicture}[x=1cm,y=1cm,
  vertex/.style={circle,fill=black,inner sep=1.45pt},
  ordinary/.style={line width=.65pt}]
  \coordinate (u) at (0,0);
  \coordinate (w) at (4.9,0);
  \coordinate (a) at (1.25,1.05);
  \coordinate (b) at (1.25,.35);
  \coordinate (c) at (1.25,-.38);
  \coordinate (d) at (1.25,-1.12);
  \coordinate (p) at (3.73,-.43);
  \coordinate (q) at (3.73,-1.08);
  \draw[densely dashed,thick,rounded corners=5pt]
    (.96,.13)--(.96,1.28)--(5.18,.2)--(5.18,-.22)--cycle;
  \draw[dotted,very thick,rounded corners=7pt]
    (.96,-.18) rectangle (3.97,-1.36);
  \foreach \x in {a,b,c,d}{\draw[ordinary] (u)--(\x);}
  \foreach \x in {p,q}{\draw[ordinary] (w)--(\x);}
  \foreach \x in {u,w,a,b,c,d,p,q}{\node[vertex] at (\x) {};}
  \node[left] at (u) {\(u\)};
  \node[right,xshift=9pt] at (w) {\(w\)};
  \node[right,yshift=-3pt] at (a) {\(a\)};
  \node[right] at (b) {\(b\)};
  \node[right] at (c) {\(c\)};
  \node[right] at (d) {\(d\)};
  \node[left] at (p) {\(p\)};
  \node[left] at (q) {\(q\)};
  \node at (2.7,1.07) {\(F\)};
  \node at (2.5,-1.59) {\(K\)};
\end{tikzpicture}
\caption{The hypergraph \(Q\). Solid segments are ordinary edges.
The dashed and dotted boundaries indicate \(F\) and \(K\), respectively.}
\label{fig:obstruction}
\end{figure}

Let \(Q\) have vertex set \(\{u,w,a,b,c,d,p,q\}\), ordinary
edges \(ua,ub,uc,ud,wp,wq\), and the two further edges
\(F=\{a,b,w\}\) and \(K=\{c,d,p,q\}\).
Figure~\ref{fig:obstruction} distinguishes the ordinary
edges from the two nonordinary edges used in the proof.

\begin{lemma}\label{lem:eight-example}
The hypergraph \(Q\) belongs to \(\Ctwo\setminus\Htwo\).
\end{lemma}

\begin{proof}
In a proper two-coloring of the ordinary edges, the vertices
\(a,b,c,d\) have the color opposite to \(u\), and \(p,q\)
have the color opposite to \(w\). The edge \(F\) is proper
precisely when \(u,w\) have the same color, whereas \(K\)
is proper precisely when they have different colors.
Thus \(Q\) is not two-colorable.
To establish criticality, we next check the effect of every
possible edge deletion, considering the two nonordinary edges first. If \(F\) is deleted, give
\(u,w\) different colors and all leaves the color opposite
to their ordinary neighbor. If \(K\) is deleted, use the
same rule with \(u,w\) of the same color. If \(ua\) or
\(ub\) is deleted, give \(u\) color \(0\), \(w\) color \(1\),
and every leaf the color opposite to its ordinary neighbor.
Then \(K\) is proper. Change the color of the freed leaf
\(a\) or \(b\) to \(0\), which also makes \(F\) proper.
For a deletion from \(\{uc,ud,wp,wq\}\), give \(u,w\)
color \(0\), all leaves color \(1\), and change the
color of the freed leaf to \(0\). The edge \(F\) remains
proper, and \(K\) now has both colors. Each leaf whose color is changed
belongs only to its deleted ordinary edge and the
hyperedge made proper by that color change. Thus none of the
other retained edges is affected by the recoloring. Hence every
edge deletion is two-colorable. As \(Q\) has no
isolated vertices, Corollary~\ref{cor:edge-deletion-criterion}
shows that it is \(3\)-critical.

For the connectivity bound, the vertices \(u,w\) have
degrees four and three, respectively, whereas every
other vertex has degree two. The singleton cuts at the degree-two vertices reduce the
remaining connectivity verification to the pair of higher-degree vertices.
The cut
\(\{u,a,b,c,d\}\mid\{w,p,q\}\) has just \(F,K\) as
crossing edges and separates \(u,w\). Every other pair
contains a degree-two vertex and is separated by its
singleton cut of size two. Therefore \(\lambda(Q)\le2\).

To exclude a representation as a Haj\'os join, we now determine
the complete list of mixed separators. First let \(z\notin\{u,w\}\). After
shrinking \(z\), each ordinary star remains connected,
and each of the two hyperedges still meets both stars.
Deleting either hyperedge therefore leaves the hypergraph
connected. Deleting an ordinary edge incident with \(z\)
has no further effect after shrinking \(z\), since that
edge would reduce to a singleton. Deleting any other
ordinary edge detaches one leaf from its star. Its
hyperedge joins it to the other star, while the other
hyperedge joins the two remaining stars. Thus
\((Q-e)\div z\) is connected for every edge \(e\).

If \(z=u\), the ordinary edges at \(u\) disappear.
The star \(wp,wq\) remains, with \(F\) attaching \(a,b\)
and \(K\) attaching \(c,d\). Deleting \(wp\) or \(wq\)
leaves its freed leaf attached through \(K\) to the
other leaf and hence to \(w\). Deleting \(F\) isolates
\(a,b\), and deleting \(K\) isolates \(c,d\).
If \(z=w\), the star at \(u\) remains, \(F\) becomes
the ordinary edge \(ab\), and \(K\) attaches \(p,q\).
Deleting \(F\) preserves connectedness. Deleting a
surviving ordinary edge leaves its leaf attached through
\(ab\) or \(K\), as appropriate. Deleting \(K\) isolates
\(p,q\). The preceding cases show that the complete list of mixed separators is
\[
\begin{array}{c|c|c}
 (z,e)&
 \text{isolated vertices of }(Q-e)\div z&
 \varphi(e)\text{ in a coloring of }Q-e\\ \hline
 (u,F)&a,b&1-\varphi(u)\\
 (u,K)&c,d&1-\varphi(u)\\
 (w,K)&p,q&1-\varphi(w)
\end{array}
\]
where \(\varphi(e)\) denotes the common color on \(e\).
In each row, all unlisted vertices form one connected
component. Thus each mixed separator in the table has exactly
three components: the two displayed isolated vertices
and the remaining connected component. Each admits
three unordered partitions into two nonempty groups;
Lemma~\ref{lem:inverse} applies to all of them through
the same invariant \(\varepsilon_{z,e}\).
The color relations follow explicitly from
the surviving hyperedge. In \(Q-F\), the edge \(K\)
forces \(u,w\) to have different colors, so \(a,b,w\)
all have the color opposite to \(u\). In \(Q-K\),
the edge \(F\) forces \(u,w\) to have the same color,
so \(c,d,p,q\) all have the color opposite to both.

The hypergraph \(Q\) is not a hyperwheel, since it has
two nonordinary edges. If it belonged to \(\Htwo\),
a generating expression would end in a Haj\'os join
of two members of \(\Htwo\subseteq\Ctwo\).
The identified vertex and new edge would form a mixed
separator whose colors agree, by Lemma~\ref{lem:inverse}.
Every mixed separator in the table has
\(\varepsilon_{z,e}=1\). The inverse criterion therefore
excludes every possible grouping of its components,
and hence every possible last Haj\'os join. It follows
that \(Q\notin\Htwo\).
\end{proof}

\Needspace{7\baselineskip}
\subsection{\texorpdfstring{Haj\'os generation and small orders}{Haj\'os generation and small orders}}\label{subsec:hajos-starting}\label{subsec:small-orders}

We now compare the two generation rules more closely.
The next lemma places every starting hypergraph in the
family generated by hyperwheels and Haj\'os joins.
We then identify a split with a singleton trace as a
Haj\'os join of the same two factors. Together, these facts will
restrict the possible splits in a recursive decomposition of a smallest obstruction.

\begin{lemma}\label{lem:hajos-bases}
All starting hypergraphs belong to \(\Htwo\).
\end{lemma}

\begin{proof}
Joining \(C_{2r-1}\) and a triangle along ordinary edges,
with the identified vertex omitted from the new edge,
gives \(C_{2r+1}\). Starting with \(C_3\), this generates
every ordinary odd cycle. In a one-node tree construction,
the node has private vertices and label \(1\), so it is
unmarked. Its node edge consists of all private vertices,
each joined to the apex by a spoke. The degree condition
gives at least three private vertices. Thus this
construction is a hyperwheel.

For the remaining tree constructions, induct on their
order. Since the tree has an edge, its proper labeling
has a node \(t\) of label \(0\). This node has no private
vertices, and \(|B_t|\ge2\) implies \(|N_T(t)|\ge2\).
Partition \(N_T(t)\) into two nonempty sets \(I_1,I_2\).
Delete the node edge \(e_t\). For each \(i\), let \(J_i\)
be the union of the branches \(G_{s\to t}\) with
\(s\in I_i\), retaining their boundary vertices and
the apex \(v\). Put
\[
 D_i=\{v\}\cup\{p_{ts}:s\in I_i\},\qquad H_i=J_i+D_i.
\]
The two \(J_i\) share only \(v\). Each \(H_i\) has
strictly fewer vertices than the original hypergraph,
because the omitted group contains a boundary vertex.
Every vertex other than \(v\) still has degree two.
At a boundary vertex, \(D_i\) replaces the deleted
edge \(e_t\); all other incidences are unchanged.
The edge \(D_i\) is new. An old edge containing a
boundary vertex \(p_{ts}\) is the node edge \(e_s\),
which does not contain \(v\) because \(\sigma(s)=1\),
whereas \(v\in D_i\).

Fix \(v\) to have color \(0\). Since every neighbor \(s\)
of \(t\) has label \(1\), equation
\eqref{eq:base-branch-colors} forces every intact
boundary vertex \(p_{ts}\) to have color \(0\).
Consequently \(H_i\) is not two-colorable. Deleting
\(D_i\) permits a coloring obtained by combining the
branch colorings, which agree at \(v\).
If an old edge \(a\in E(J_i)\) is deleted, take a proper
two-coloring of the original hypergraph with \(a\)
deleted and normalize the color of \(v\) to \(0\).
All branches in the other group are intact, so their
boundary vertices have color \(0\). The original
edge \(e_t\) is retained and proper. Since \(P_t\)
is empty, some boundary vertex in the \(i\)-th group
must therefore have color \(1\). Restricting this coloring to the vertices of the factor properly colors
\(H_i-a\), including its added edge \(D_i\). There are no isolated
vertices in \(H_i\). Applying
Corollary~\ref{cor:edge-deletion-criterion} to each factor
therefore shows that both \(H_i\) are \(3\)-critical.

Each \(H_i\) has at most one vertex of degree greater
than two. Lemma~\ref{lem:bases} therefore identifies it as
an ordinary odd cycle or a labeled-tree construction. The odd
cycles are already in \(\Htwo\), and the smaller tree
constructions belong to \(\Htwo\) by induction.
Join disjoint copies of \(H_1,H_2\) at their apices,
using \(D_1,D_2\) as the deleted edges. Their union
away from \(v\) is \(B_t\). Including \(v\) in the
new edge precisely when \(t\) was marked restores
\(e_t\), and hence the original hypergraph.
\end{proof}

\begin{lemma}\label{lem:singleton}
A Toft split of two members of \(\Ctwo\) at a degree-two
vertex with a singleton trace is admissible and is a
Haj\'os join of its two factors.
\end{lemma}

\begin{proof}
A singleton is split by no vertex partition, so every
cut splits at most one of the two traces. The first
factor is \(3\)-critical with maximum local edge
connectivity at most two, and therefore each pair of
its vertices has a separating cut of size two.
This proves admissibility.
Suppose \(S_1=\{s\}\), and write the incident edges at
the split vertex as \(\{v\}\cup A_1,\{v\}\cup A_2\).
Form a Haj\'os join by identifying \(s\in V(G_1)\)
with \(v\in V(G_2)\) and deleting \(e\) and
\(\{v\}\cup A_2\). The retained incident edge becomes
\(\{s\}\cup A_1=S_1\cup A_1\).
Since \(S_1\cup S_2=e\), the set \(S_2\) equals
\(e\setminus\{s\}\) or \(e\). Omit or include the
identified vertex in the new edge according as
\(s\notin S_2\) or \(s\in S_2\). That edge is then
exactly \(S_2\cup A_2\). The vertex and edge sets
are those of the prescribed split.
\end{proof}

The small-order classification requires keeping track of
which vertices can gain degree in a split. A direct comparison
of incidences in a degree-two split gives
\begin{equation}\label{eq:split-degrees}
 d_G(x)=d_{G_1}(x)+\mathbf{1}_{\{x\in S_1\cap S_2\}}
 \quad(x\in V(G_1)),\qquad
 d_G(y)=d_{G_2}(y)
 \quad(y\in V(G_2)\setminus\{v\}).
\end{equation}
Indeed, a vertex of \(e\) loses that edge and gains one
or two inserted edges, and the traces cover \(e\).
Other vertices of the first factor lose no edge.
Each incidence at a vertex of the second factor is
retained or replaced by one incidence in an inserted edge.

Write \(H_3\) for the hypergraph whose ordinary edges
form the path \(xavby\) and whose remaining edge is
\(\{x,v,y\}\), following the notation in
\cite[Proposition~7.2 and Figure~1]{sstorch}.
It is the three-node tree construction with one
private vertex at each leaf and its middle node marked.

\begin{lemma}\label{lem:rim-and-path-traces}
Suppose an admissible split has \(|S_1|,|S_2|\ge2\).
If the selected edge is the rim of a hyperwheel, the
two traces are disjoint. If the first factor is \(H_3\)
and the selected edge is \(\{x,v,y\}\), then
\(\{S_1,S_2\}=\{\{x,v\},\{v,y\}\}\).
\end{lemma}

\begin{proof}
The traces cover the selected edge and have intersection
of size at most one, by
Corollary~\ref{cor:admissibility-test}. In particular,
the selected edge has size at least three.
For a hyperwheel with apex \(u\), suppose a rim vertex
\(r\) belongs to both traces. In any cut separating
\(r,u\), consider the side avoiding \(u\), which is
a nonempty set of rim vertices. If it contains
\(j\) vertices and is a proper subset of the rim,
the crossing edges are its \(j\) spokes and the rim,
so the cut has size \(j+1\). The whole rim gives
at least three crossing spokes. Hence the only
size-two cut separating \(r,u\) isolates \(r\).
As each trace has at least two vertices, this cut
splits both traces, contrary to admissibility.
The traces must therefore be disjoint.

For the triple of \(H_3\), the two traces have size
two and share one vertex. Suppose the shared vertex
is \(x\). A cut separating \(x,v\) necessarily splits
the triple, so a size-two cut crosses just one
ordinary edge of the path \(xavby\). That edge is
\(xa\) or \(av\), because these are the edges on the
path from \(x\) to \(v\). After deleting it, \(v,y\)
lie in the same component of the remaining ordinary
path. Thus the cut restricts to \(x\mid\{v,y\}\)
on the triple and splits both traces. This is
impossible. Interchanging \(x,y\) excludes an overlap
at \(y\), leaving only the overlap at \(v\) stated
in the lemma.
\end{proof}

For use in the classification that follows, let \(W_n\) denote the
hyperwheel on \(n\) vertices.
The five-vertex assertion in the next lemma describes the restriction of
Stiebitz, Storch, and Toft's classification of all nine
five-vertex \(3\)-critical hypergraphs
\cite[Proposition~7.2 and Figure~1]{sstorch}.
Section~7 of that paper also reports a computer classification
on six vertices. For the obstruction argument, we need the
additional connectivity restriction and explicit six-vertex
edge lists. The proof below verifies the additional connectivity restriction directly
and derives the required edge lists from the recursive
characterization established in the preceding section.

\begin{lemma}\label{lem:small-members}
Up to isomorphism, the members of \(\Ctwo\) of order at most
five are \(C_3,W_4,C_5,W_5,H_3\).
There are exactly five members of order six, given by
the following edge lists. A string of distinct vertex
labels denotes their set.
\begin{equation}\label{eq:six}
\begin{array}{c|l|l|c}
 \text{type}&\text{ordinary edges}&\text{other edges}
 &\text{vertices of degree at least three}\\ \hline
 A&01,02,03,04,05&12345&0\\
 B&01,02,13,14&125,345&1\\
 C&01,02,13,24,25&345&2\\
 D&02,03,12,14,15&2345&1,2\\
 E&02,05,13,14&125,345&1,5
\end{array}
\end{equation}
\end{lemma}

\begin{proof}
Every \(3\)-critical hypergraph has at least three
vertices. On three vertices it must be the triangle.
Indeed, if an ordinary edge is missing, give its two endpoints the same color and give
the third vertex the other color.
This properly colors every possible remaining ordinary
edge and the triple. If all three ordinary edges are
present, they form a \(3\)-critical subhypergraph, which
must be the whole hypergraph.
Suppose \(G\in\Ctwo\) has at most five vertices and at
least two vertices of degree greater than two. The
decomposition in the proof of Theorem~\ref{thm:main}
gives factors in \(\Ctwo\) whose orders sum to
\(|V(G)|+1\le6\). Each factor has at least three
vertices, so both are triangles. By
\eqref{eq:split-degrees}, the vertices of degree
greater than two in the split are exactly the vertices
in \(S_1\cap S_2\). There is at most one such vertex
by admissibility, a contradiction. Thus every member
of order at most five is a starting hypergraph.

For a tree construction with \(t\) nodes and \(p\)
private vertices, the order is \(n=t+p\), and
\[
 2t\le\sum_{s\in V(T)}|B_s|=p+2(t-1).
\]
Hence \(p\ge2\) and \(t\le n-2\). A two-node tree is
impossible. Each of its leaves needs a private vertex
to have \(|B_s|\ge2\), so both adjacent nodes would have
label \(1\). With three nodes, the tree is a path.
Both leaves need private vertices and have label \(1\);
the middle node has label \(0\) and no private vertices.
For \(n\le5\), there is exactly one private vertex at
each leaf. The middle node must be marked to make
the apex degree at least three, giving \(H_3\).
A one-node tree gives \(W_4\) or \(W_5\), and the
ordinary odd cycles in this range are \(C_3,C_5\).
This proves the first assertion. It also shows directly
that every listed hypergraph belongs to \(\Ctwo\).

Now let \(n=6\). A starting hypergraph is a tree
construction with \(t\le4\), since an ordinary odd
cycle has odd order. For \(t=4\), the tree is a
path or a star. The path has leaves in opposite
label classes, whereas both leaves need private
vertices and label \(1\). The star needs a private
vertex at each of its three leaves, but \(p=6-4=2\).
Thus neither occurs. For \(t=3\), the leaves have
one and two private vertices, in either order.
The middle node can be marked or unmarked, giving
types \(B,C\), respectively. In the displayed labeling,
\(B\) has apex \(1\), private sets \(\{0\},\{3,4\}\),
and node edges \(02,125,345\).
The hypergraph \(C\) has apex \(2\), private sets
\(\{0\},\{4,5\}\), and node edges \(01,13,345\).
The one-node construction gives type \(A\).

A member of order six that is not a starting hypergraph has at least two
vertices of degree greater than two. Its factors
have orders three and four, so they are a triangle
and \(W_4\). If the selected edge in the first factor
is ordinary, admissibility forces a singleton trace.
If it is the three-vertex rim of \(W_4\), two
nonsingleton traces would overlap, contrary to
Lemma~\ref{lem:rim-and-path-traces}. Thus
Lemma~\ref{lem:singleton} identifies this split as a
Haj\'os join of the triangle and \(W_4\).
All vertices other than the identified vertex retain
their degrees in a Haj\'os join. If \(\eta\in\{0,1\}\)
indicates whether that vertex belongs to the new edge,
its degree is the sum of the two original degrees
minus two, plus \(\eta\). Identifying the apex of
\(W_4\) gives only one vertex of degree greater than
two. Identifying a rim vertex gives two such vertices
precisely when \(\eta=1\). Up to automorphisms of the
factors, the chosen edge of \(W_4\) is then its rim
or the spoke at that rim vertex. These are the only
two possibilities. Joining the wheel with apex \(1\)
and rim \(245\) to the triangle \(023\), identifying
\(2\), deleting \(245,23\), and including \(2\) in
the new edge gives type \(D\). Joining the wheel
with apex \(1\) and rim \(345\) to the triangle
\(025\), identifying \(5\), deleting \(15,25\),
and including \(5\) gives type \(E\).
Both belong to \(\Ctwo\) by Lemma~\ref{lem:hajos-input}.
Finally, the edge sizes and degree multisets distinguish
the five displayed types from one another, completing the classification.
\end{proof}

\begin{proof}[Proof of the lower bound in Theorem~\ref{thm:hajos}]
We prove that every \(G\in\Ctwo\) of order at most
seven belongs to \(\Htwo\), by induction on its order.
If \(G\) is a starting hypergraph, use
Lemma~\ref{lem:hajos-bases}. Otherwise it has at least
two vertices of degree greater than two. Choose its
decomposition from the proof of Theorem~\ref{thm:main}.
The smaller factors \(G_1,G_2\) are in \(\Ctwo\),
so both belong to \(\Htwo\) by induction.
A singleton trace gives a Haj\'os join by
Lemma~\ref{lem:singleton}, and hence puts \(G\)
in \(\Htwo\).
Assume both traces have size at least two. Admissibility
gives \(|S_1\cap S_2|\le1\), so \(|e|\ge3\).
The order relation
\(|V(G_1)|+|V(G_2)|=|V(G)|+1\), together with
\(|V(G_2)|\ge3\), gives \(|V(G_1)|\le5\).
By Lemma~\ref{lem:small-members}, the selected edge
is the rim of \(W_4\) or \(W_5\), or the triple of
\(H_3\).

For \(W_4\), the rim cannot be covered by two disjoint
nonsingleton traces, contradicting
Lemma~\ref{lem:rim-and-path-traces}. For \(W_5\),
that lemma shows that the traces form a disjoint
\(2+2\) partition of its rim. Here \(|V(G_1)|=5\), so
\(|V(G_2)|=3\) and \(G_2\) is a triangle.
Equation~\eqref{eq:split-degrees} shows that the
split has just the original apex of \(W_5\) as
a vertex of degree greater than two, a contradiction.
For \(H_3\), the traces are \(\{x,v\},\{v,y\}\)
by Lemma~\ref{lem:rim-and-path-traces}. Again
\(G_2\) is a triangle. Only \(v\) gains a degree
on the first side, and it was already the sole
vertex of degree greater than two in \(H_3\).
The surviving vertices of the triangle have degree
two. This again contradicts the choice of \(G\).
The preceding cases exhaust the possibilities and complete the induction.
Lemma~\ref{lem:eight-example} supplies an obstruction
of order eight, so the minimum order is eight.
\end{proof}

\subsection{Uniqueness at order eight}

The remaining finite analysis concerns a six-vertex
first factor and a triangle as second factor. We carry out
this analysis in a separate lemma, giving the full list of possibilities
and the exclusion of each alternative before using the classification
in the proof of uniqueness.

\begin{lemma}\label{lem:six-trace-classification}
Let \(G_1\in\Ctwo\) have six vertices and let
\(G_2\) be a triangle. Form an admissible Toft split
\(G\) with \(|S_1|,|S_2|\ge2\), and suppose that
\(G\) has at least two vertices of degree greater
than two. Up to an automorphism of the labeled first
factor and simultaneous interchange of the traces
and the two surviving vertices \(6,7\) of \(G_2\),
the possibilities are exactly the three rows of
\eqref{eq:tracecases}. For each row, the displayed
mixed separator gives a Haj\'os decomposition with
factor orders \((4,5)\), \((4,5)\), or \((3,6)\),
respectively. In particular, every such \(G\)
belongs to \(\Htwo\).
\end{lemma}

\begin{proof}
By Lemma~\ref{lem:small-members}, the first factor is
one of the five types in \eqref{eq:six}; this list is
exhaustive by the preceding proof. Use its vertex labels
for \(G_1\), and label the two surviving vertices of
\(G_2\) as \(6,7\). The split retains \(67\)
and inserts \(S_1\cup\{6\},S_2\cup\{7\}\).
Equation~\eqref{eq:split-degrees} shows that its
vertices of degree greater than two are exactly
those already of degree greater than two in \(G_1\),
together with the vertices in \(S_1\cap S_2\).
We examine every edge of size at least three in
each of the five types. Since the traces cover that
edge and each has size at least two, no ordinary edge
can occur. Corollary~\ref{cor:admissibility-test} gives
\(|S_1\cap S_2|\le1\). Thus a selected triple has two
two-element traces with one common vertex; a selected
four-element edge has either disjoint two-element
traces or traces of sizes two and three with one
common vertex. Type \(A\), whose rim has size five,
is treated first below by applying the hyperwheel trace lemma.
These possibilities exhaust the trace sizes.

For the case analysis, we use the following description of size-two cuts.
For a connected hypergraph \(J\), call an edge
\(h\in E(J)\) a \emph{separating edge} if \(J-h\)
is disconnected. This condition allows several vertices
of \(h\) to remain in the same component after its deletion.
The hypergraph \(G_1-e\) is connected, since every
nontrivial cut of a \(3\)-critical hypergraph has size
at least two. If a cut of \(G_1\) has size two and
splits \(e\), it crosses exactly one edge \(h\) of
\(G_1-e\). Deleting \(h\) leaves no edge of \(G_1-e\)
crossing that partition. Thus \(h\) must be a separating
edge of \(G_1-e\), and every component of
\((G_1-e)-h\) lies wholly on one side of the cut.
In each of the cases considered below, deleting such
an edge leaves exactly two components. Their vertex
sets therefore determine the partition, up to
interchanging its sides. We check whether this
partition separates the specified pair of vertices
and splits \(e\); when it does, its crossing edges
in \(G_1\) are precisely \(e,h\). Cuts that leave
\(e\) intact are considered separately.

For type \(A\), Lemma~\ref{lem:rim-and-path-traces}
shows that the two traces must be disjoint. The split then still
has only one vertex of degree greater than two,
which contradicts the choice of \(G\).
For types \(B,C\), there is one such vertex
in \(G_1\), so the traces must overlap in a
degree-two vertex to produce a second one.
Admissibility will exclude each of these possible
overlaps in the two types under consideration.

In type \(B\), first take \(e=125\). The only
possible new vertices of degree greater than two
are the overlap vertices \(2,5\). Both lie in
\(e\), as does the existing vertex \(1\) of
degree greater than two. A cut separating either
overlap vertex from \(1\) must split \(e\).
The separating edges of \(G_1-e\) are \(01,02,345\).
The other two edges \(13,14\) each have an
alternate connection between their endpoints
through \(345\), so neither is a separating edge.
Deleting \(01\) or \(02\) partitions \(e\) as
\(2\mid15\), and deleting \(345\) partitions
it as \(5\mid12\). Thus every size-two cut
separating the overlap from \(1\) isolates
that overlap on its side of \(e\), splitting
both traces. This violates admissibility.
For \(B\) with \(e=345\), an overlap \(r\) is
one of \(3,4,5\). If a cut separating \(r,1\)
left \(e\) intact, the edges \(13,14,125\)
would all cross, so it would have size at
least three. Otherwise its one crossing edge
in \(G_1-e\) must be a separating edge. The edges
\(01,02\) each have an alternate connection
through \(125\), so the separating edges are exactly
\(13,14,125\). Their deletions isolate
\(3,4,5\), respectively, and give
the partitions \(3\mid45\), \(4\mid35\),
\(5\mid34\), respectively. Only the partition
isolating \(r\) separates it from \(1\), so
every relevant cut again splits both traces.

For type \(C\), the selected edge is \(e=345\),
and the existing vertex of degree greater than
two is \(2\). An overlap \(r\) is one of
\(3,4,5\). Leaving \(e\) intact on the side
opposite to \(2\) crosses \(24,25\) and
at least one edge of the path \(2013\).
If \(e\) is split, \(G_1-e\) is the ordinary
tree with edges \(01,02,13,24,25\).
Its unique crossing edge lies on the path
from \(2\) to \(r\). For \(r=3\), the
possible edges are \(02,01,13\), and all
give \(3\mid45\) on \(e\). For \(r=4,5\),
the only possibilities are \(24,25\),
giving \(4\mid35\), \(5\mid34\), respectively.
Thus every size-two cut separating \(r,2\)
splits both traces, excluding type \(C\).

For type \(D\), the selected edge is \(e=2345\).
Every size-two cut separating \(1,2\) splits \(e\).
If \(e\) were intact, it would be opposite to \(1\),
and \(12,14,15\) would all cross. Once \(e\) is
split, the only crossing edge of the ordinary tree
\(G_1-e\) is \(12\), giving \(23\mid45\) on \(e\).
If the traces are disjoint, each has size two.
To avoid splitting both, they must be \(23,45\).
If the traces overlap, their sizes are two and three.
The larger trace necessarily crosses \(23\mid45\),
so the smaller trace must be \(23\) or \(45\).
With smaller trace \(23\), an overlap at \(3\)
is impossible. Every size-two cut separating \(2,3\)
splits \(e\) and crosses either \(02\) or \(03\),
the edges of the path \(203\). Both choices give
\(3\mid245\) on \(e\), splitting both traces.
Thus the overlap must be \(2\), giving \(23,245\).
With smaller trace \(45\), let \(r\in\{4,5\}\)
be the overlap. A size-two cut separating \(1,r\)
cannot leave \(e\) intact, because \(12,14,15\)
would all cross. If \(e\) is split, its one
ordinary crossing edge is \(1r\), giving
\(r\mid(e\setminus\{r\})\) on \(e\).
Both traces are split, excluding this possibility.
There are exactly the two stated cases for type \(D\).

For type \(E\), each selected edge is a triple,
so the traces have size two and share one vertex.
First take \(e=125\). A cut separating the vertices
\(1,5\) of degree greater than two splits \(e\).
The separating edges of \(G_1-e\) are \(02,05,345\).
The remaining edges \(13,14\) each have an
alternate connection through \(345\).
Only the separating edge \(345\) separates \(1,5\). Its cut gives
\(1\mid25\) on \(e\), excluding overlap at \(1\).
An overlap at \(2\) is excluded by separating \(2,5\).
The only possible crossing separating edges are \(02,05\),
both giving \(2\mid15\) and splitting both traces.
The remaining overlap is \(5\), with traces \(15,25\).
Now take \(e=345\) in type \(E\). The separating edges
of \(G_1-e\) are \(13,14,125\). The edges
\(02,05\) each have an alternate connection
through \(125\), so neither is a separating edge.
If a cut separating
\(1,5\) left \(e\) intact, \(13,14,125\) would all
cross. If it splits \(e\), the only crossing edge
of \(G_1-e\) that separates \(1,5\) is \(125\),
giving \(5\mid34\) on \(e\). This excludes overlap
at \(5\). For an overlap \(r\in\{3,4\}\), consider
a cut separating \(r,1\). Leaving \(e\) intact
again crosses \(13,14,125\). Splitting \(e\)
leaves only the ordinary edge \(1r\) as a possible
crossing edge of \(G_1-e\), since \(r\) is a leaf
there. It gives \(r\mid(e\setminus\{r\})\),
splitting both traces. Thus this choice of selected edge
does not yield any admissible case under consideration.

After these exclusions, the complete list is
\eqref{eq:tracecases}, up to automorphisms of the first
factor and interchanging \(6,7\) together with the two
traces. The last column records the two component vertex
sets obtained after deleting \(f\) and shrinking \(z\).
We will use these partitions to recover a Haj\'os join
in each remaining case.
\begin{equation}\label{eq:tracecases}
\begin{array}{c|c|c|c|c|c}
 \text{type}&e&S_1&S_2&(z,f)&
 \text{components of }(G-f)\div z\\ \hline
 D&2345&23&45&(1,67)&457\mid0236\\
 D&2345&23&245&(2,2457)&145\mid0367\\
 E&125&15&25&(5,257)&02\mid13467
\end{array}
\end{equation}
We next verify that all three remaining cases are admissible. In either type \(D\)
row, the set \(W\) from
Corollary~\ref{cor:admissibility-test} is
\(\{1,2\}\), and the cut \(023\mid145\)
has crossing edges \(12,2345\) and splits at most
one trace. For the type \(E\) row, \(W=\{1,5\}\),
and the cut \(134\mid025\), with crossing edges
\(125,345\), has the same property.

To apply the inverse criterion, we verify that the
last column gives exactly the components for each
mixed separator in the table. In the first row, deleting \(67\)
and shrinking \(1\) leaves edges \(457\) on
one component and \(02,03,236\) on the other.
In the second row, deleting \(2457\) and
shrinking \(2\) leaves \(14,15\) on one
component and the path edges \(03,36,67\)
on the other. In the last row, deleting
\(257\) and shrinking \(5\) leaves \(02\)
on one component and \(13,14,34,16,67\)
on the other. These edge lists show that each displayed partition
is precisely the partition into connected components.

In the last two rows, \(z\in f\), so the common
color of \(f\) in every proper two-coloring
of \(G-f\) equals the color of \(z\).
For the first row, \(14,15\) give \(4,5\)
the color opposite to \(1\). The retained
edge \(457\) then forces \(7\) to have
the color of \(1\). Since \(G\) is critical,
\(67\) is monochromatic in every coloring
of \(G-67\), so its common color also
equals the color of \(1\).
Lemma~\ref{lem:inverse} applies in all three
rows and gives factors of orders \((4,5)\),
\((4,5)\), and \((3,6)\), respectively.
Both factors in each row are \(3\)-critical with
maximum local edge connectivity at most two. Their
orders are at most six, so the lower bound already
proved in Theorem~\ref{thm:hajos} places them in
\(\Htwo\). Their Haj\'os join is \(G\), and hence
\(G\in\Htwo\) in every row.
\end{proof}

\begin{proof}[Proof of the uniqueness assertion of Theorem~\ref{thm:hajos}]
Let \(G\in\Ctwo\setminus\Htwo\) have eight vertices.
It has at least two vertices of degree greater than
two, by Lemmas~\ref{lem:bases} and~\ref{lem:hajos-bases}.
Use the decomposition from the proof of
Theorem~\ref{thm:main}, with factors \(G_1,G_2\).
Both factors have smaller order and hence belong
to \(\Htwo\) by the lower bound just proved.
Neither trace is a singleton, by
Lemma~\ref{lem:singleton}. Therefore \(|e|\ge3\)
and the factor orders sum to nine.
The first factor cannot have order three.
If it has order four, it is \(W_4\), whose rim
is excluded by Lemma~\ref{lem:rim-and-path-traces}.
The only remaining ordered pairs of factor orders are
\((5,4)\) and \((6,3)\). We consider these two cases separately,
retaining the order of the factors throughout because the traces
are subsets of an edge in the first factor.

Suppose first that their orders are \((5,4)\).
Then \(G_2=W_4\). Write its apex as \(w\), its
degree-two split vertex as \(z\), and its other
rim vertices as \(p,q\). The sets at the split
vertex are \(A_1=\{w\}\), from the spoke \(wz\),
and \(A_2=\{p,q\}\), from the rim.
The retained edges on this side are \(wp,wq\).
Among the five-vertex members of \(\Ctwo\), only
\(W_5,H_3\) have an edge of size at least three.
If \(G_1=W_5\), its traces form a disjoint \(2+2\)
partition of the rim. Write its apex as \(u\),
and its traces as \(\{a,b\},\{c,d\}\), where
the first is paired with \(A_1\). The inserted
edges are \(\{a,b,w\}\) and \(\{c,d,p,q\}\).
Together with the retained spokes, these are
exactly the edges of \(Q\).

For the other possibility, take \(G_1=H_3\), with traces
\(\{v,x\},\{v,y\}\). Interchanging the ends of
the path if necessary, pair \(\{v,x\}\) with
\(A_1=\{w\}\). The split has ordinary edges
\(xa,av,vb,by,wp,wq\) and further edges
\[
 f=\{v,x,w\},\qquad \{v,y,p,q\}.
\]
After deleting \(f\) and shrinking \(v\), the
components have vertex sets \(\{x,a\}\) and
\(\{b,y,w,p,q\}\). The first contains the edge
\(xa\), and the second is connected through
\(by,wp,wq,\{y,p,q\}\). Since \(v\in f\),
the common color of \(f\) in \(G-f\) equals
the color of \(v\). Lemma~\ref{lem:inverse}
therefore gives critical factors of orders three
and six, both in \(\Ctwo\). By the lower bound already established, both belong to
\(\Htwo\). Their join would therefore imply
\(G\in\Htwo\), a contradiction.

We now turn to the remaining pair \((6,3)\).
Lemma~\ref{lem:six-trace-classification} applies: the
split is admissible, both traces have size at least
two, and \(G\) has at least two vertices of degree
greater than two. That lemma gives \(G\in\Htwo\),
again a contradiction.
Thus every eight-vertex member of
\(\Ctwo\setminus\Htwo\) arises from the
\((5,4)\) case producing \(Q\).
Lemma~\ref{lem:eight-example} verifies that
\(Q\) itself belongs to this difference,
which proves uniqueness up to isomorphism.
\end{proof}

\begin{corollary}\label{cor:two-traces}
The class \(\Ctwo\) is generated from hyperwheels by
both forms of the Haj\'os join and by admissible Toft
splits at degree-two vertices in which both traces
have size at least two.
\end{corollary}

\begin{proof}
Hyperwheels belong to \(\Ctwo\), and both types of
operation preserve membership in \(\Ctwo\), by the
closure assertion of Lemma~\ref{lem:hajos-input} and
Theorem~\ref{thm:main}. This establishes one inclusion in the claimed generation rule.
For the reverse inclusion, take a decomposition of a member
of \(\Ctwo\) supplied by Theorem~\ref{thm:main}.
Replace each starting hypergraph by its Haj\'os
construction from Lemma~\ref{lem:hajos-bases}.
At each internal node with a singleton trace,
Lemma~\ref{lem:singleton} replaces the split
by a Haj\'os join of the same two factors.
Every remaining split has both traces of size
at least two. Induction up the finite decomposition
tree gives a construction using precisely the
operations in the statement.
\end{proof}

\section{Related results and consequences}\label{sec:related}
\label{subsec:cut-functions}\label{subsec:edge-counts}

The transfer formula can also be expressed through an
edge splitting function. In the setting of
Lemma~\ref{lem:transfer}, the prescribed traces define
the following function on the subsets of the selected
edge \(e\),
\[
 \theta(A)=k-\max\left\{
       |\{i:S_i\subseteq A\}|,
       |\{i:S_i\subseteq e\setminus A\}|
                   \right\},
 \qquad A\subseteq e.
\]
Equation~\eqref{eq:transfer} says that the minimum size
of a cut extending \(U\mid V(G_1)\setminus U\) is
\(|\cut_{G_1-e}(U)|+\theta(U\cap e)\).
The two counts inside the maximum correspond to the two possible placements
of the vertex representing the contracted second factor side relative to the cut.
The function is symmetric under complementation in
\(e\), with \(\theta(\varnothing)=\theta(e)=0\).
If \(\varnothing\ne A\subsetneq e\), neither \(A\)
nor its complement can contain all the traces,
because \(\bigcup_{i=1}^k S_i=e\). Both counts are
therefore at most \(k-1\), giving
\(1\le\theta(A)\le k\). This is an edge splitting
function in the framework of Veldt, Benson, and
Kleinberg \cite{vbk}.

Local connectivity preservation also appears in the
splitting-off theorem of B\'erczi, Chandrasekaran,
Kir\'aly, and Kulkarni, which gives a complete
splitting-off preserving local edge connectivities
in weighted hypergraphs \cite{bckk}.
Cho and Oh construct terminal mimicking networks
preserving minimum multicut values up to a
prescribed threshold \cite{co}.
For the prescribed Toft split considered here,
Corollary~\ref{cor:split-connectivity} gives a
necessary and sufficient condition in terms of
size-\(k\) cuts of the first factor and the given
traces, when both factors have maximum local edge
connectivity at most \(k\).
Lemma~\ref{lem:transfer} also gives the exact
preservation of local edge connectivities on
\(V(G_2)\setminus\{v\}\) and the value \(k\)
for pairs on opposite sides. For \(k=2\), the
criterion is the admissibility test of
Corollary~\ref{cor:admissibility-test}.

The recursive characterization also determines the exact
edge count of every member of \(\Ctwo\). We first derive
this count from the starting hypergraphs and the effect of the split
operation. We then relate the resulting equality to the general
lower bound on the number of edges in a critical hypergraph.

\begin{corollary}\label{cor:edges}
Every \(H\in\Ctwo\) satisfies \(|E(H)|=|V(H)|\).
\end{corollary}

\begin{proof}
Every ordinary odd cycle has equally many edges
and vertices. In a tree construction with \(t\)
nodes and \(p\) private vertices, the tree has
\(t-1\) edges and hence \(t-1\) assigned vertices.
Together with the apex and private vertices,
this gives \(1+(t-1)+p=t+p\) vertices.
There is one node edge for each tree node
and one spoke for each private vertex.
These \(t+p\) edges are all distinct, as
verified in Section~\ref{sec:bases}.
Thus every starting hypergraph has equal
edge and vertex counts.

Consider a degree-two Toft split of two
vertex-disjoint critical factors \(G_1,G_2\).
It deletes the selected edge \(e\) of \(G_1\)
and exactly the two edges incident with
the split vertex \(v\) of \(G_2\).
The two inserted edges \(S_1\cup A_1\) and
\(S_2\cup A_2\) are distinct. Indeed,
their intersections with
\(V(G_2)\setminus\{v\}\) are \(A_1,A_2\),
which are distinct because the original
edges \(\{v\}\cup A_1,\{v\}\cup A_2\)
are distinct. Both inserted edges meet
both factor sides, since all \(S_i,A_i\)
are nonempty. They therefore coincide
with none of the retained edges, each
of which lies wholly on one side.
The retained edge families from the
two factors are also disjoint.
Consequently the split \(G\) satisfies
\begin{align*}
 |V(G)|&=|V(G_1)|+|V(G_2)|-1,\\
 |E(G)|&=(|E(G_1)|-1)+(|E(G_2)|-2)+2
        =|E(G_1)|+|E(G_2)|-1.
\end{align*}
Equality of edge and vertex counts for
both factors therefore implies equality
for the split. Theorem~\ref{thm:main}
expresses every member of \(\Ctwo\)
by a finite sequence of admissible
splits from the starting hypergraphs.
Induction over this construction proves
the corollary.
\end{proof}

The edge-count equality also follows from a sparsity
property of blocks, without using a recursive generating
expression for the hypergraph. The next lemma gives this property for
every connected hypergraph with no separating vertex
and maximum local edge connectivity at most two.
Its strict inequality applies to every proper nonempty
edge family, while the whole edge set satisfies a
nonstrict inequality.

\begin{lemma}\label{lem:block-sparsity}
Let \(B\) be a connected hypergraph with at least
two vertices, no separating vertex, and
\(\lambda(B)\le2\). For an edge family \(F\), put
\(V(F)=\bigcup_{e\in F}e\). Then
\begin{equation}\label{eq:block-sparsity}
 |F|<|V(F)|\quad
   (\varnothing\ne F\subsetneq E(B)),\qquad
 |E(B)|\le |V(B)|.
\end{equation}
In particular, every nonempty proper vertex
set \(U\subsetneq V(B)\) contains at most
\(|U|-1\) edges of \(B\).
\end{lemma}

\begin{proof}
Suppose there is a nonempty edge family with
at least as many edges as vertices in its
support, and choose one, \(F\), minimal under
inclusion. It has at least two edges, since
each edge has at least two vertices. Put
\(m=|F|\), \(t=|V(F)|\). For every \(a\in F\),
minimality gives
\[
 m-1<|V(F\setminus\{a\})|\le t\le m.
\]
The quantities are integers, so \(m=t\),
and every proper nonempty subfamily of \(F\)
has strictly fewer edges than vertices.
Let \(K=(V(F),F)\). For a partition
\(V(F)=X\cup Y\) into nonempty parts, the
family of edges contained in either part
is a proper subfamily of \(F\). If it is
nonempty, its size is less than its support
size, and hence at most \(|X|-1\) or
\(|Y|-1\), respectively. The same bounds
hold for an empty family. Since
\(|F|=|X|+|Y|\), at least two edges cross
this partition. Thus \(K\) is connected
and every nontrivial cut has size at
least two.

To apply the block lemma to this subhypergraph, we must also rule out
the presence of a separating vertex. If \(z\) separated \(K\), partition the
components of \(K\div z\) into two nonempty
groups. Each edge of \(K\) is supported
on one group together with \(z\), so this
partitions \(F\) into two nonempty proper
families \(F_1,F_2\). Each group contributes
an edge because every vertex of \(K\)
lies in an edge, and edges have size
at least two. Their supports cover \(V(F)\)
and meet in at most \(z\). Minimality gives
\[
 |F|=|F_1|+|F_2|
 \le |V(F_1)|+|V(F_2)|-2
 \le |V(F)|-1,
\]
a contradiction. Hence \(K\) has no
separating vertex.
Lemma~\ref{lem:embedded-block}, applied inside
\(B\), shows that \(K\) is a block of \(B\).
The hypergraph \(B\) itself is connected
and has no separating vertex, so maximality
forces \(K=B\). Consequently any
inclusion-minimal family with at least
as many edges as support vertices is
the whole edge set of \(B\) and has
equal edge and vertex counts. There
can therefore be no such proper family,
and \(|E(B)|\le|V(B)|\). If no such family
\(F\) exists at the outset, then every nonempty
edge family already satisfies the strict
inequality. This proves
\eqref{eq:block-sparsity} in both cases.

Finally, let \(U\) be a nonempty proper
vertex set and let \(F_U\) consist of
the edges contained in \(U\). Connectedness
and the absence of isolated vertices in
\(B\) imply \(F_U\subsetneq E(B)\).
If \(F_U\ne\varnothing\), then
\(|F_U|<|V(F_U)|\le|U|\).
If \(F_U=\varnothing\), its size is
zero and is again at most \(|U|-1\).
This proves the last assertion.
\end{proof}

Combining Lemma~\ref{lem:block-sparsity} with the
following classical lower bound for critical hypergraphs
gives a second proof of Corollary~\ref{cor:edges}.

\begin{theorem}[{\cite{lovasz1970}}]\label{thm:critical-edge-input}
Every \(3\)-critical hypergraph \(H\) satisfies
\(|E(H)|\ge|V(H)|\). For every integer \(n\ge3\),
there is a \(3\)-critical hypergraph with exactly
\(n\) vertices and \(n\) edges.
\end{theorem}

This result is due to Lov\'asz \cite{lovasz1970};
see also \cite{lovasz}. Later proofs were given by Woodall
\cite{woodall} and Seymour \cite{seymour}.
Its attribution and extremal form are recorded in
\cite[Section~8, p.~445]{sstorch}.
For \(H\in\Ctwo\), criticality gives at least
three vertices, connectedness, and the
absence of a separating vertex. These facts allow
Lemma~\ref{lem:block-sparsity} to be applied,
giving the reverse inequality.
Together with Theorem~\ref{thm:critical-edge-input},
the sparsity argument therefore also yields the
edge-count equality.
A hypergraph is called \emph{square} if its numbers
of edges and vertices are equal. Thus every member
of \(\Ctwo\) is square.
Abbott and Hare \cite{ah} studied square
\(3\)-critical hypergraphs.
The square condition alone, however, does not force the maximum
local edge connectivity of a critical hypergraph to be at most two.
For this comparison, we use \(H_5\) from the classification
of Stiebitz, Storch, and Toft
\cite[Proposition~7.2 and Figure~1]{sstorch}.
Its criticality and edge count are already known.
The additional calculation below determines its maximum
local edge connectivity.

\begin{lemma}[{\cite[Proposition~7.2 and Figure~1]{sstorch}}]
\label{lem:known-h5}
Let \(H_5\) have vertex set \(\{x,y,z,p,q\}\), ordinary
edges \(xz,zy,pq\), and further edges
\(F=\{x,y,p\}\) and \(K=\{x,y,q\}\).
Then \(H_5\) is \(3\)-critical and
\(|E(H_5)|=|V(H_5)|=5\).
\end{lemma}

With the vertices of \(H_5\) relabeled as above,
its two degree-three vertices determine the maximum
local edge connectivity; every other pair contains a
degree-two vertex. The next lemma uses this distinction
to compare the square condition with the additional
requirement \(\lambda\le2\).

\begin{lemma}\label{lem:square-example}
The hypergraph \(H_5\) of Lemma~\ref{lem:known-h5}
satisfies \(\lambda(H_5)=3\).
\end{lemma}

\begin{proof}
Use the vertex and edge labels of
Lemma~\ref{lem:known-h5}. We first examine the pair
of degree-three vertices, and then use singleton
cuts to bound every other pair.
Every cut separating \(x,y\) crosses
both \(F,K\), since both contain
\(x,y\). It also crosses an edge of
the ordinary path \(xzy\), because
a path joining opposite sides of
a vertex partition must have a
crossing edge. These are three
distinct edges, so
\(\lambda_{H_5}(x,y)\ge3\).
The singleton cut at \(x\) has size
\(d_{H_5}(x)=3\), giving equality.
Finally \(z,p,q\) all have degree
two. Every pair other than
\(\{x,y\}\) therefore admits a
singleton cut of size two.
Hence \(\lambda(H_5)=3\), as claimed.
\end{proof}

\Needspace{7\baselineskip}
\section*{Declaration of generative AI and AI-assisted technologies}
The author used a large language model during the preparation of this
manuscript for assistance with mathematical exploration, checking
arguments, organization and language editing. The author is responsible
for the mathematical statements, proofs, references and final text.

\end{document}